\documentclass[12pt]{article}
\usepackage{amsfonts}
\usepackage{amssymb}
\usepackage{amsmath}
\usepackage{eufrak}
\usepackage[dvips]{graphicx}
\usepackage{latexsym}

\usepackage[all]{xy}

\newcommand{\C}{{\EuFrak C}}

\newcommand{\G}{{\cal G}}

\newcommand{\bi}{\begin{itemize}}
\newcommand{\ei}{\end{itemize}}

\newtheorem{theorem}{Theorem}[section]
\newtheorem{lemma}[theorem]{Lemma}
\newtheorem{fact}[theorem]{Fact}

\newtheorem{corollary}[theorem]{Corollary}
\newtheorem{proposition}[theorem]{Proposition}
\newtheorem{definition}[theorem]{Definition}
\newtheorem{remark}[theorem]{Remark}
\newtheorem{conjecture}[theorem]{Conjecture}
\newtheorem{question}[theorem]{Question}

\newtheorem{example}[theorem]{Example}

\def\defn{\begin{definition}\upshape}
\def\edefn{\end{definition}}

\newtheorem{main theorem}{Theorem}
\newtheorem{main theorem repeated}{Theorem}
\newtheorem{main conjecture}[main theorem]{Conjecture}
\newtheorem{main conjecture repeated}[main theorem repeated]{Conjecture}

\title{Borel equivalence relations and  Lascar strong types}
\author{Krzysztof Krupi\'nski\footnote{Research supported by the Polish Government grant N N201 545938}  \hspace{0.1pt}, Anand Pillay\footnote{Supported by EPSRC grant EP/I002294/1} and S\l awomir Solecki\footnote{Supported by  NSF grant DMS-1001623}}
\date{April 16, 2012}

\begin{document}
\maketitle
\begin{abstract}
The ``space" of {\em Lascar strong types}, on some sort and relative to a given complete theory $T$, is in general {\em not} a compact Hausdorff topological space. We have at least three (modest) aims in this paper. The first is to show that  spaces of Lascar strong types, as well as other related spaces and objects such as the Lascar group $Gal_{L}(T)$ of $T$, have  well-defined {\em Borel cardinalities} (in the sense of the theory of complexity of Borel equivalence relations). The second is to compute the Borel cardinalities of the known examples as well as of some new examples that we give. The third is to explore notions of {\em definable} map, embedding, and isomorphism, between these and related quotient objects. We also make some conjectures, the main one being roughly  ``smooth iff trivial"
The possibility of a descriptive set-theoretic account of the complexity of spaces of Lascar strong types was touched on in the paper \cite{CLPZ}, where the first example of a ``non $G$-compact theory" was given. The motivation for writing this paper is partly the discovery of new examples via definable groups, in \cite{CoPi1}, \cite{CoPi2} and the generalizations in \cite{GiKr}.

\end{abstract}
\footnotetext{2010 Mathematics Subject Classification: 03C45, 03E15}
\footnotetext{Key words and phrases: Lascar strong types, Borel reducibility}

\section[\mbox{}]{Introduction}

Strong types of one form or another play an important role in the study of first order theories $T$. For example classical (or Shelah) strong types form an obstruction to ``uniqueness of nonforking extensions", and when $T$ is stable they form the only obstruction (the finite equivalence relation theorem). Likewise for so-called KP strong types and the independence theorem in simple theories. The most general form of strong types are what have been called {\em Lascar strong types}. In contradistinction to Shelah strong types and KP strong types, the collection of Lascar strong types (in a given sort) does not in general have the structure of a (Hausdorff) topological space (although it has a natural but maybe degenerate ``quasi-compact" topology). So one issue is how to view a space of Lascar strong types as a mathematical object.
In \cite{CLPZ} it was suggested that Lascar strong types could or should be looked at from the point of view of descriptive set theory, in particular the general theory of quotients of Polish spaces by Borel equivalence relations. We pursue this theme in the current paper, trying at least to give robust definitions and analyzing some examples. A brief remark (for the set-theorists) is that there {\em are} some fundamental compact spaces attached to a first order theory; namely the type spaces. The ``space" of countable models of a theory $T$ is in a sense a special case, being a $G_{\delta}$-subspace of  $S_{\omega}(T)$ the space of complete types of $T$ in variables $x_{1},x_{2},\dots$. In any case, type spaces and their quotients provide many interesting objects attached to a first order theory $T$, and this paper is in a sense concerned with how to understand ``bad" quotients. \\

We briefly recall these various strong types. We fix a complete, possibly many-sorted, theory $T$, a monster model $\C$, and a ``small" set $A$ of parameters from $\C$. Let us fix also a sort, or even $A$-definable set $S$, and variables $x, y, \dots$ of that sort. Elements $a,b$ of sort $S$ are said to have the same {\em Shelah strong type} over $A$ if $E(a,b)$ whenever $E(x,y)$ is an $A$-definable equivalence relation with only finitely many classes. The equivalence relation of having the same Shelah strong type over $A$ is denoted  $E_{S,A,Sh}$. The collection of Shelah strong types over $A$, i.e. $E_{S,A,Sh}$ equivalence classes, is a Stone space (compact, totally disconnected) and corresponds to the space of complete types in sort $S$ over $acl^{eq}(A)$.  We say that elements $a,b$ from sort $S$ have the same {\em KP strong type over $A$}, if $E(a,b)$ whenever $E(x,y)$ is an equivalence relation given by maybe infinitely many formulas over $A$ and with at most $2^{|A|+ |T|}$ classes. We denote the  equivalence relation by $E_{S,A,KP}$. The collection of KP strong types over $A$ has the natural structure of a compact Hausdorff space (but not necessarily totally disconnected). The nomenclature ``KP" stands for Kim-Pillay and was introduced by Hrushovski in \cite{Hrushovski}. We stick with the notation as it has become accepted now. Finally we say that elements $a,b$ of sort $S$ have the same {\em Lascar strong type} if  $E(a,b)$ whenever $E$ is an $Aut({\C}/A)$-invariant equivalence relation with at most $2^{|A|+|T|}$-many classes. We denote the equivalence relation by $E_{S,A,L}$. As remarked above, the collection of Lascar strong types, i.e. $E_{S,A,L}$-classes, does not have the natural structure of a Hausdorff topological space. It is clear that  $E_{S,A,KP}$ refines $E_{S,A,Sh}$ and that $E_{S,A,L}$ refines $E_{S,A,KP}$. \\

We could restrict each of these equivalence relations to a type-definable over $A$ set $X$ to get  $E_{X,A,Sh}$ etc. (And it is a fact from \cite{Lascar-Pillay} for example that the restriction of $E_{S,A,KP}$ to $X$ coincides with the finest type-definable over $A$ equivalence relation with at most $2^{|A| + |T|}$ classes). 
We can also consider these equivalence relations on some possibly infinite product of sorts $\prod_{i}S_{i}$, or on some $*$-definable set over $A$ (type-definable over $A$ set of possibly infinite tuples). One of our main objects of study will be the restriction of $E_{S,A,L}$ to an $E_{S,A,KP}$ class $X$.

We often drop $S$ when it is clear. And when $A = \emptyset$, it is omitted. 

The general idea is to investigate properties of $E_{L}$ in terms of Borel equivalence 
relations on type spaces. We assume $T$ is countable so that the relevant type spaces (over countable models) are Polish spaces.

\begin{definition}
Let $E$ and $F$ be Borel equivalence relations on Polish spaces $X$ and $Y$, respectively. We say that $E$ is Borel reducible to $F$ if there exists a Borel reduction of $E$ into $F$, i.e. a Borel function $f \colon X \to Y$ such that for all $x,y \in X$
$$xEy \iff f(x)Ff(y).$$
If $E$ is Borel reducible to $F$, we write $E\leq_{B}F$. We say that $E$ and $F$ are Borel equivalent (or Borel bi-reducible), symbolically $E \sim_B F$, if $E \leq_B F$ and $F \leq_B E$. We will write $E \lneq_B F$ if $E \leq_B F$ and $F \nleq_B E$.

$E$ is said to be smooth if it is Borel reducible to the equality relation on a Polish space (equivalently on the Cantor set), namely there is a Borel function $f$ from $X$ to a Polish space $Y$ such that for all $x,y\in X$, $f(x) = f(y)$ iff $xEy$.

Sometimes we identify a  Borel equivalence relation $E$ on $X$ with the ``space" $X/E$, and speak of $X/E \leq_B Y/F$. 
\end{definition}

We will point out in Section 3 how to consider the set of $E_{L}$-classes as the quotient of a type space over a countable model by a Borel (in fact $K_{\sigma}$) equivalence relation, which in the sense of Borel equivalence above does not depend on the model chosen. 

Typically we will be investigating the restriction of $E_{L}$ to a given KP type over $\emptyset$, $X$ (namely $E_{KP}$-equivalence class) which we note as $E_L \! \upharpoonright \! X$. Recall that in this case Newelski proved that $E_L \! \upharpoonright \! X$ is either trivial (i.e. $X$ is itself an $E_{L}$-class) or has exactly continuum many classes \cite[Corollary 1.8(2)]{Ne}. From the Borel reducibility perspective, this means that if $E_L$ is non-trivial on $X$, then equality on the Cantor set (denoted by $\Delta_{2^{\mathbb N}}$) is Borel reducible to $E_L \! \upharpoonright \! X$. By the Silver dichotomy \cite{Si}, $\Delta_{2^{\mathbb N}}$ is the simplest (in the sense of Borel reducibility) uncountable Borel equivalence relation. So, a natural question arises whether $E_L\! \upharpoonright \! X$ can be Borel equivalent with $\Delta_{2^{\mathbb N}}$ (and so smooth), or whether it must be essentially more complicated (non-smooth).  Recall that Harrington-Kechris-Louveau dichotomy \cite{HKL} tells us that each Borel equivalence relation $E$ is either smooth  or the relation $E_0$ is Borel reducible (even embeddable) to $E$, where $E_0$ is the relation on $2^{\mathbb N}$ defined by
$$xE_0y \iff (\exists n) (\forall m \geq n)(x(m)=y(m)).$$ Namely, $E_{0}$ is the simplest non-smooth Borel equivalence relation.
%
The following conjecture is aimed at strengthening Newelski's theorem mentioned above:

\begin{main conjecture}\label{main conjecture}
Let $X$ be an $E_{KP}$-class. Then  $E_L \! \upharpoonright \! X$ is either trivial or non-smooth (so, in the latter case, $E_0$ is Borel reducible to $E_L \! \upharpoonright \! X$).
\end{main conjecture}

Another justification for the above conjecture is the fact that if  $E_L \! \upharpoonright \! X$ is non-trivial, then all classes of $E_L \! \upharpoonright \! X$ are not $G_\delta$ subsets \cite[page 167]{Ne}. Recall that we have $E_0\leq_B E$ for every orbit equivalence relation $E$ induced by a continuous action of a Polish group on a Polish space whose all classes are non-$G_\delta$; see \cite[Corollary 3.4.6]{BeKe}.

On the other hand, we will see shortly that $E_L \!\upharpoonright \! X$ is always a $K_\sigma$ equivalence relation (countable union of compacts).  There is a most complicated $K_{\sigma}$ equivalence relation. One realization of such an equivalence relation is $ l^\infty$ \cite[Theorem 6.6.1]{Ka}, namely the  relation on ${\mathbb R}^{\mathbb N}$ defined by
$$x  l^\infty y \iff \sup_n |x(n)-y(n)| < \infty.$$
%

In this paper, we give examples where (up to Borel equivalence) $E_L \! \upharpoonright \! X$ is $E_{0}$ and where it is $l^{\infty}$. We do not know examples with other Borel cardinalities. Nevertheless:
\begin{main conjecture} Any non-smooth $K_{\sigma}$ equivalence relation can be represented by some $E_L \! \upharpoonright \! X$ (in some theory $T$).
\end{main conjecture}

We will also mention the {\em Lascar group}  $Gal_{L}(T)$ of a theory $T$. This is precisely the group of permutations of all Lascar strong types (even of countable tuples) over $\emptyset$ induced by $Aut(\C)$. In the case of KP strong types, we have $Gal_{KP}(T)$ which is naturally a compact Hausdorff (separable) group. 
There is a canonical surjective homomorphism $\pi:Gal_{L}(T) \to Gal_{KP}(T)$. We let $Gal_{0}(T)$ denote $ker(\pi)$.  $Gal_{0}(T)$ (as $Gal_{L}(T)$ itself) can again be viewed in a robust manner as the quotient of a Polish space by a $K_{\sigma}$ equivalence relation as we describe in Section 5.  $T$ is said to be $G$-compact when $Gal_{0}(T)$ is trivial. $Gal_{0}(T)$ is a fundamental invariant of the bi-interpretability type of $T$, and a special case of Conjecture 1 is:

\begin{main conjecture} $Gal_{0}(T)$ is either trivial or non-smooth.
\end{main conjecture}


%

Many of the examples studied in this paper come from definable groups. Given a $\emptyset$-definable group $G$ (in $\C\models T$), and some small set $A$ of parameters, we have the ``connected components" $G_{A}^{0}$, $G_{A}^{00}$ and $G_{A}^{000}$. 

\begin{itemize}
\item  $G_{A}^{0}$ is the intersection of all $A$-definable subgroups of $G$ of finite index. 
\item  $G_{A}^{00}$ is the smallest type-definable over $A$ subgroup of $G$ of ``bounded" index (which amounts to index at most $2^{|A|+|T|}$). 
\item $G_{A}^{000}$ is the smallest subgroup of $G$ of bounded (as above) index which is also $Aut(\C/A)$-invariant. 
\end{itemize}

These are normal subgroups of $G$ and moreover, $G^{000}_A \leq  G^{00}_A \leq  G^{0}_A \leq G$.  When $A = \emptyset$, we will omit it. But actually in many cases (e.g. when $T$ is $NIP$), these connected components do not depend on the choice of $A$ anyway.

The quotients $G/G^{0}$, $G/G^{00}$ and $G/G^{000}$ are group analogues of the spaces of Shelah strong types, KP strong types and Lascar strong types. $G/G^{0}$ is naturally a profinite group, $G/G^{00}$ is naturally a compact Hausdorff group and we again want to treat $G/G^{000}$  (as well as $G^{00}/G^{000}$) as a descriptive set theoretic object. We can and will make Conjecture 1 for $G^{00}/G^{000}$. 
In fact, as is well-known, there is a close connection with the various strong types. Namely,  if we add a new sort $S$ and a regular action of $G$ on $S$, then $S$ is itself a complete type and the group quotients above correspond to $E_{S,Sh}$, $E_{S,KP}$ and $E_{S,L}$, respectively. In particular, $G^{00}/G^{000}$ will be Borel equivalent to $X/E_{S,L}$ as we point out below (for $X$ any $E_{S,KP}$ class).


In Section 1, we give descriptive set-theoretic and model-theoretic preliminaries. In Section 2, we explain how to interpret our various quotient structures in a robust manner as descriptive set-theoretic objects. In Section 3, we 
compute the Borel cardinalities of the existing examples as well as some new examples witnessing non $G$-compactness, with $G^{00}/G^{000}$ playing a prominent role. 
Section 4 focuses on Borel cardinalities in a certain axiomatic framework from \cite{GiKr} (generalizing examples from \cite{CoPi1}), and confirms Conjecture \ref{main conjecture} in the classes of concrete examples obtained in \cite{GiKr}. In Section 5, we point out that the ``Lascar group" also has a robust descriptive set-theoretic character and we make some observations and conjectures. Finally, in Section 6 we explore notions of ``definable" map, embedding, isomorphism, between our various quotient structures, and extend results of Thomas 
\cite{Thomas} on Borel reductions which are not continuous reductions.



\section{Preliminaries}

First, we recall relevant notions and facts from descriptive set theory. 

A Polish space is a separable, completely metrizable topological space. A Borel equivalence relation on (a Polish space) $X$ is an equivalence relation on $X$ which is a Borel subset of $X \times X$. A function $f \colon X \to Y$ is said to be Borel if the preimage of every Borel subset of $Y$ is a Borel subset of $X$.

For a Polish space $X$, let $\Delta_X$ be the (equivalence relation of) equality  on $X$. Then, we have the following relations between some basic Borel equivalence relations:
$$\Delta_1 \lneq_{B} \Delta_2 \lneq_B \dots \lneq_B \Delta_{\mathbb N} \lneq_B \Delta_{2^{\mathbb N}} \lneq_B E_0.$$ 
Moreover, we have the following two fundamental dichotomies \cite{Si, HKL}.

\begin{fact}[Silver dichotomy]\label{Silver dichotomy}
For every Borel equivalence relation $E$ either $E \leq_B \Delta_{\mathbb N}$ or $\Delta_{2^{\mathbb N}} \leq_B E$. 
\end{fact}

\begin{fact}[Harrington-Kechris-Louveau dichotomy]\label{Harrington-Kechris-Louveau dichotomy}
For every Borel equivalence relation $E$ either $E \leq_B \Delta_{2^{\mathbb N}}$ (in which case $E$ is smooth) or $E_0 \leq_B E$. 
\end{fact}

The linearity of $\leq_B$ breaks drastically above $E_0$, but we will not discuss it here. The interested reader is referred to \cite{Ka}.

Recall that a $K_\sigma$ subset of a topological space is a countable union of compact subsets of this space. For example, all the above basic relations as well as the relation $l^\infty$ defined in the introduction are  $K_\sigma$.  
The next fact says that $l^\infty$ is universal for $K_\sigma$ relations \cite[Theorem 6.6.1]{Ka}.

\begin{fact}\label{universality of l^infty}
Each $K_\sigma$-equivalence relation is Borel reducible to $l^\infty$.
\end{fact}

We will need the following theorem \cite[Theorem 3.4.5]{BeKe}.

\begin{fact}\label{dense orbit}
Let $X$ be a Polish space and $G$ a Polish group acting on $X$ by homeomorphisms. Let $E_{a}$ denote the orbit equivalence relation on $X$. Assume there exists a dense orbit under this action and that $E_a$ is a meager subset of $X\times X$. Then $E_0 \leq_B E_a$.  
\end{fact}

Let $X$ be a Polish space. By $F(X)$ we denote the collection of all closed subsets of $X$. We endow $F(X)$ with the $\sigma$-algebra generated by the sets
$$\{ F \in F(X) : F \cap U \ne \emptyset \},$$
where $U$ varies over open subsets of $X$. The set $F(X)$ equipped with this $\sigma$-algebra is called the Effros-Borel space of $F(X)$. It turns out that there exists a Polish topology on $F(X)$ whose Borel sets form the above $\sigma$-algebra \cite[Theorem 12.6]{Ke}. If $X$ is additionally compact, this property is fulfilled by the so-called Vietoris topology, that is the topology whose basis is formed by the sets 
$$\{ F \in F(X) : F \cap K =\emptyset \wedge F \cap U_1 \ne \emptyset \wedge \dots \wedge F \cap U_n \ne \emptyset \}$$
with $K$ ranging over closed subsets and $U_1, \dots, U_n$ over open subsets of $X$ \cite[Exercise 12.7]{Ke}.

We will need the following fact \cite[Theorem 12.13]{Ke}.

\begin{fact}\label{function d}
Let $X$ be a Polish space.  Then there exists a Borel mapping $d \colon F(X) \to X$ such that $d(F) \in F$ for every non-empty $F \in F(X)$.
\end{fact}

\begin{definition}
Let $E$ be an equivalence relation on a set $X$. For any $A \subseteq X$, by $[A]_E$ (the so-called saturation of $A$) we denote the set of all elements of $X$ whose $E$-class meets $A$. A selector for $E$ is a map $s \colon X \to X$ such that whenever $xEy$, then $s(y)=s(x) \in [x]_E$. A transversal for $E$ is a set $T\subseteq X$ meeting every equivalence class of $E$ in exactly one point.
\end{definition}



We will need the following fact, whose part (ii) will be a typical method of computing Borel complexity. Point (i) of the fact can be found in 
\cite[Exercise 24.20 and p. 363]{Ke}; point (ii) follows immediately from point (i). 
 
\begin{fact}\label{Borel section}
(i) Each surjective, continuous function $f \colon X \to Y$ between compact, Polish spaces $X$ and $Y$ has a Borel section, i.e. a Borel function $g \colon Y \to X$ for which $f\circ g = id_{Y}$.
\newline
(ii) Hence, if $X$, $Y$ are compact Polish spaces, $E, F$ are Borel equivalence relations on $X$, $Y$, respectively, and $f\colon X \to Y$ is a continuous surjective function such that $x_{1}Ex_{2}$ iff $f(x_{1})Ff(x_{2})$ for all $x_{1},x_{2}\in X$, then $E \sim_{B} F$.
\end{fact}
%


Now, we turn to the model-theoretic preliminaries. We fix a complete countable theory $T$ and work in a so-called monster model $\C$ of $T$. We let $x,y$ range over a given sort or even a $*$-definable over $\emptyset$ set. 
Let $\Theta(x,y)$ be the $\emptyset$-type-definable relation expressing that the $2$-element sequence $(x,y)$ extends to an infinite indiscernible sequence.  The following are well-know facts: 

\begin{fact}\label{Theta}
$E_L$ (on the given sort) is the transitive closure of $\Theta$. In particular, if the language is countable, then $E_L$ is a union of countably many $\emptyset$-type-definable sets.
\end{fact}

\begin{fact}\label{G^{000}}
Let $G$ be a $\emptyset$-definable group or even a group $*$-definable over $\emptyset$. Then
$G^{000}$ is the subgroup of $G$ generated by $\{ gh^{-1} : g,h \in G, \; E_L(g,h)\}$.  In particular, if the language is countable, then by Fact \ref{Theta}, $G^{000}$ is a union of countably many $\emptyset$-type-definable sets.
\end{fact}

Now, we discuss in more detail the relationship between model-theoretic connected components and Lascar strong types (e.g. see \cite[Section 3]{GiNe}).

 Let $M$ be a first order structure, and $G$ a group definable in $M$; without loss of generality we assume that it is $\emptyset$-definable. Consider the 2-sorted structure $N=(M,X,\cdot)$ described by the conditions:
\begin{enumerate}
\item[i)] $M$ and $X$ are disjoint sorts,
\item[ii)] $\cdot \colon G \times X \to X$ is a regular action of $G$ on $X$,
\item[iii)] $M$ is equipped with its original structure.
\end{enumerate} 
Consider a monster model $N^*=(M^*,X^*,\cdot) \succ N$ of $Th(N)$, where $M^*\succ M$ is a monster model of $Th(M)$. Let $G^*$ be the interpretation of $G$ in $M^*$, and ${G^*}^{00}$ and ${G^*}^{000}$ the connected components of $G^*$ computed in $M^*$ (equivalently in $N^*$). 
Note that all elements of $X^*$ have the same type over $\emptyset$ in $N^*$.
The following fact relates connected components of $G^*$ with Lascar strong types and Kim-Pillay types on the sort $X^*$ of $N^*$ \cite[Lemma 3.7]{GiNe}.
\begin{fact}\label{regular action}
Let $x \in X^*$. Then:\\
(i) The $E_{L}$ equivalence class $[x]_{E_L}$ of $x$ in $N^*$ coincides with its ${G^*}^{000}$-orbit.\\
(ii) The $E_{KP}$ equivalence class  $[x]_{E_{KP}}$ of $x$ in $N^*$ coincides with its ${G^*}^{00}$-orbit. 
\end{fact}
This shows that any example with ${G^*}^{000} \ne {G^*}^{00}$ yields a new example of a non-$G$-compact theory, namely $Th(N^*)$. In the next section, we will use the above fact to show that the relation on ${G^*}^{00}$ of lying in the same orbit modulo ${G^*}^{000}$, considered on a certain space of types, is Borel equivalent to $E_L \! \upharpoonright \! [x]_{E_{KP}}$ for any $x \in X^*$.

It had been an open question whether there exists a group $G$ such that  $G^{000} \ne G^{00}$. This was answered in the affirmative in \cite{CoPi1}. One of the examples found there is (a saturated model of the theory of) the universal cover of $SL_2(\mathbb R)$. Subsequently, \cite{GiKr} describes a certain general situation in which an extension $\widetilde{G}$ of a group $G$ by an abelian group $A$ by means of a 2-cocycle with finite image satisfies $\widetilde{G}^{000} \ne \widetilde{G}^{00}$, yielding new classes of concrete examples with this property (e.g. some central extensions of $SL_2(k)$ or, more generally, of higher dimensional symplectic groups $Sp_{2n}(k)$, for ordered fields $k$). More details on this will be given in Section 4, where we analyze the Borel cardinality of $E_L$, verifying Conjecture \ref{main conjecture} in some general situations from \cite{GiKr}.

\section{$\pmb{E_{L}}$ as a Borel equivalence relation}

The goal of this section is to explain how to treat both $E_L$ and the relation of lying in the same coset modulo $G^{000}$ (where $G$ is a $\emptyset$-definable group) as Borel (even $K_\sigma$) equivalence relations in reasonably canonical ways, so as to formalize the discussion and conjectures from the introduction. \\


As before we work in a monster model $\C$ of a first order theory $T$ in a countable language.

We consider the relation $E_L$ on a sort $X$ of $\C$ (or even on infinite tuples of elements from some sorts). We will write $S(M)$ having in mind $S_X(M)$. An important fact (see for example \cite[Fact 1.9]{CLPZ}) is that if (for a model $M\prec \C$) $tp(a/M) = tp(b/M)$ then $E_{L}(a,b)$.  In fact, to be more precise, $tp(a/M) = tp(b/M)$ implies there is $c$ such that $\Theta(a,c)$ and $\Theta(c,b)$; and if $\Theta(a,b)$, then $tp(a/M) = tp(b/M)$ for some model $M$. 

\begin{definition}\label{Lascar as Borel}
Let $M\prec \C$ be a countable submodel. We define the binary relation $E^M_L$ on $S(M)$ by 
$$pE^M_L q \iff (\exists a \models p) (\exists b \models q)(a E_L b).$$  
\end{definition}

\begin{remark}\label{very basic}
(i) $pE^M_L q \iff (\forall a \models p) (\forall b \models q)(a E_L b)$.\\
(ii) $E^M_L$ is a $K_\sigma$ equivalence relation on $S(M)$. \\
(iii) The map taking $a\in X$ to $tp(a/M)$ induces a bijection between $X/E_{L}$ and $S_{X}(M)/E^M_L$. 
\end{remark}
{\em Proof.} (i) $(\Leftarrow)$ is obvious. For $(\Rightarrow)$, choose $a_0 \models p$ and $b_0 \models q$ such that $a_0 E_L b_0$, and consider any $a \models p$ and $b \models q$. Since $tp(a/M)=tp(a_0/M)$ and $tp(b/M)=tp(b_0/M)$, we get (as remarked above) $aE_L a_0$ and $b E_L b_0$, and so $a E_L b$.\\
(ii) follows from (i) and Fact \ref{Theta}.\\
(iii) As remarked above, the quotient map  $X \to X/E_{L}$ factors through the map taking $a\in X$ to $tp(a/M)$, and the resulting equivalence relation on $S(M)$ is clearly $E^M_L$. 
 \hfill $\blacksquare$\\
 
There are some slightly subtle issues around Remark 2.2. The equivalence relation $E_{L}$ on tuples corresponds, on the face of it, to a $K_{\sigma}$ subset of the space $S_{X\times X}(M)$, and the latter is {\em not} the product space $S_{X}(M)\times S_{X}(M)$. But because the $E_{L}$-class of a tuple depends only on the type of that tuple over $M$, in fact $E_{L}$ corresponds to a $K_{\sigma}$ equivalence relation on the space $S_{X}(M)$.
 
In any case, Definition 2.1 and Remark 2.2 describe  how  the space $X/E_{L}$ is viewed as a Polish space quotiented by a Borel equivalence relation. The following  proposition shows that the resulting ``Borel complexity" does not depend on the countable model chosen. 

\begin{proposition}\label{independence of a model}
Let $M$ and $N$ be countable, elementary substructures of $\C$. Then $E^M_L \sim_B E^N_L$. 
\end{proposition}
{\em Proof.} It is enough to prove the proposition under the assumption that $M\prec N$. Let $f \colon S(N) \to S(M)$ be the restriction map.

$f$ is continuous and surjective. 
It is immediate from Definition \ref{Lascar as Borel} and Corollary \ref{very basic}(i) that $E^N_L(p,q)$ iff $E^M_L(f(p), f(q))$  (as realizations $a$ of $p$ and $b$ of $q$ such that $E_{L}(a,b)$ are also realizations of $f(p)$ and $f(q)$, respectively).
By Fact \ref{Borel section}(ii), $E^N_L\sim_{B} E^M_{L}$.\hfill $\blacksquare$\\

Analogous observations hold after restriction to a single KP type. We give some details. 

\begin{definition}
Let $M\prec \C$ be a countable submodel. We define a binary relation $E^M_{KP}$ on $S(M)$ by 
$$pE^M_{KP} q \iff (\exists a \models p) (\exists b \models q)(a E_{KP} b).$$  
\end{definition}

\begin{remark}
(i) $pE^M_{KP} q \iff (\forall a \models p) (\forall b \models q)(a E_{KP} b)$.\\
(ii) $E^M_{KP}$ is a closed equivalence relation on $S(M)$. In particular, the equivalence classes of $E^M_{KP}$ are closed and so compact subsets of $S(M)$.\\
(iii) The map taking $x\in X$ to $tp(a/M) \in S_{X}(M)$ induces a bijection between 
$X/E_{KP}$ and $S_{X}(M)/E^{M}_{KP}$ (and, in fact, a homeomorphism, where $X/E_{KP}$ is given the ``logic topology"). 
\end{remark}

Let $Y$ now denote an $E_{KP}$-class, namely the $E_{KP}$-class of $a$ for some $a\in X$. (So $Y$ is no longer type-definable over $\emptyset$, but is type-definable over $bdd^{heq}(\emptyset)$ for those who know what that means.) Let $[tp(a/M)]_{E^M_{KP}}$ be what it says (the $E^M_{KP}$-class of $tp(a/M)$). Then, as in 2.2 (iii), the map taking $b\in Y$ to $tp(b/M)$ induces a bijection between $Y/E_{L}$ and 
$[tp(a/M)]_{E^M_{KP}}/E^{M}_{L}$, and we identify in this way $Y/E_{L}$ with the quotient of a compact Polish space by a Borel equivalence relation. 

A similar argument to the proof of Proposition \ref{independence of a model} yields the following observation, saying that the resulting Borel complexity of $Y/E_{L}$ is well-defined (does not depend on the countable model $M$ chosen). 

\begin{proposition}\label{independence of a model 2}
Let $M$ and $N$ be any countable, elementary substructures of $\C$. Then, for any $a$, $E^M_L\! \! \upharpoonright [tp(a/M)]_{E^M_{KP}} \sim_B E^N_L \!\!\upharpoonright [tp(a/N)]_{E^N_{KP}}$. 
\end{proposition}

Since Propositions \ref{independence of a model} and  \ref{independence of a model 2} tell us that the Borel cardinality of both $E^M_L$ and $E^M_L \! \! \upharpoonright [tp(a/M)]_{E^M_{KP}}$ does not depend of the choice of $M$, we usually skip the letter $M$ in $E^M_L$ and $E^M_{KP}$, having in mind that formally $E_L$ and $E_{KP}$  denote $E^M_L$ and $E^M_{KP}$ for some countable model $M \prec \C$. We hope it will be clear from the context when $E_L$ denotes $E^M_L$ for some $M$ and when it has its original meaning.

It was mentioned in the introduction that Newelski proved in \cite{Ne} that whenever $[a]_{L} \ne [a]_{KP}$, then $[a]_{KP}$ is refined into  continuum many equivalence classes of $E_L$. This means that $[tp(a/M)]_{E_{KP}}$ is refined into continuum many equivalence classes of $E_L$.
Using Fact \ref{Silver dichotomy}, Newelski's theorem can be reformulated as follows.
\begin{fact}\label{Newelski}
If $[a]_{L} \ne [a]_{KP}$, then $\Delta_{2^{\mathbb N}} \leq_B E_L\!\! \upharpoonright [tp(a/M)]_{E_{KP}}$.
\end{fact}
More formally, Conjecture \ref{main conjecture} should be written in the following way.
\begin{main conjecture repeated}
If $[a]_{L} \ne [a]_{KP}$, then $E_0 \leq_B E_L\!\! \upharpoonright [tp(a/M)]_{E_{KP}}$.
\end{main conjecture repeated}

We now look at the analogous notions for definable groups. Assume $G$ is group $\emptyset$-definable in $\C$.
We first explain  how $G^{00}/G^{000}$ can be seen as the quotient of a type space by a Borel equivalence relation. 

\begin{definition}
Let $M\prec \C$ be a countable submodel. We define the binary relation $E^M_{G}$ on $S_{G^{00}}(M)$ by 
$$pE^M_G q \iff (\exists a \models p) (\exists b \models q)(ab^{-1} \in G^{000})$$  
\end{definition}

\begin{remark}\label{equivalence relation}
(i) $pE^M_G q \iff (\forall a \models p) (\forall b \models q)(ab^{-1} \in G^{000})$.\\
(ii) $E^M_G$ is a $K_\sigma$-equivalence relation on $S_{G^{00}}(M)$. 
\newline
(iii) The map taking $g\in G^{00}$ to $tp(g/M)$ induces a bijection between $G^{00}/G^{000}$ and $S_{G^{00}}(M)/E^M_G$. 
\end{remark}
{\em Proof.}  (i) $(\Leftarrow)$ is obvious. For $(\Rightarrow)$, choose $a_0 \models p$ and $b_0 \models q$ such that $a_0 b_0^{-1}\in G^{000}$, and consider any $a \models p$ and $b \models q$. Since $tp(a/M)=tp(a_0/M)$ and $tp(b/M)=tp(b_0/M)$, we get $aE_L a_0$ and $b E_L b_0$. This implies that $a a_0^{-1} \in G^{000}$ and $b b_0^{-1} \in G^{000}$ (because the relation of lying in the same coset modulo $G^{000}$ is bounded and $\emptyset$-invariant, and so it is coarser than $E_L$). 
Thus, $ab^{-1}=aa_0^{-1}a_0b_0^{-1}b_0b^{-1} \in G^{000}$.\\
(ii) follows from (i) and Fact \ref{G^{000}}. \\
(iii)  As earlier with $E_{L}$.
\hfill $\blacksquare$\\

A similar argument to the proof of Proposition \ref{independence of a model} yields the following observation.

\begin{proposition}\label{independence of a model 3}
Let $M$ and $N$ be any countable, elementary substructures of $\C$. Then $E^M_G \sim_B E^N_G$. 
\end{proposition}
As in the case of $E_L$, because of the above proposition, we will usually write $E_G$ instead of $E^M_G$. \\

Everything we have said above for $\emptyset$-definable groups $G$ also holds for $*$-definable (over $\emptyset$) groups. 
Anyway, what we have done so far shows that some of the basic objects we are considering in this paper:  $G^{00}/G^{000}$, and  $X/E_{L}$ where $X$ is an $E_{KP}$-class, can be assigned a ``Borel cardinality (or complexity)" in a coherent manner. 
We will discuss similar issues for the ``Lascar group" in Section 5. 

\begin{proposition}\label{E_G sim E_L}
Consider the structure $N^{*}$ defined  before Fact 1.10 and we use notation from there. Let $x\in X$, and let $Z$ be its $E_{KP}$-class. Then $(G^{*})^{00}/(G^{*})^{000} \sim_{B} Z/E_{L}$. 
\end{proposition}
{\em Proof.} 
We can assume that $N=(M,X,\cdot)$ is countable (recall that $M \prec M^*$ and $x \in X$). 
In this proof, we will be very precise with the notation, i.e. the relation $E^M_G$ is considered on $S_{{G^*}^{00}}(M)$ (working in $M^*$), and $E^N_L$ is considered on $S_{X^*}(N)$ (working in $N^*$).

Define the function $f \colon S_{{G^*}^{00}}(M) \to S_{X^*}(N)$ by
$$f(tp(g/M)) = tp(g\cdot x/N).$$
Notice that this function is well-defined. For this, take any $g_1, g_2 \in {G^*}^{00}$ satisfying the same type over $M$. Then there exists $f \in Aut(M^*/M)$ mapping $g_1$ to $g_2$. But $f$ gives rise to an $\overline{f} \in Aut(N^*/N)$ which is defined by the conditions $\overline{f}\! \upharpoonright \! M^* =f$ and $\overline{f}(g\cdot x)=f(g)\cdot x$ for $g \in G^*$. Therefore, $\overline{f}(g_1 \cdot x)=f(g_1)\cdot x=g_2 \cdot x$, and we conclude that $tp(g_1 \cdot x/N)=tp(g_2 \cdot x/N)$. 

We check that $f$ is a continuous surjection from $S_{{G^*}^{00}}(M)$ to $[tp(x/N)]_{E^N_{KP}}$ and is a Borel reduction of $E^M_G$  to $E^N_L \!\! \upharpoonright [tp(x/N)]_{E^N_{KP}}$. The fact that the range of $f$ equals $[tp(x/N)]_{E^N_{KP}}$ follows from Fact \ref{regular action}(ii). The continuity of $f$ is clear. Finally, the fact that $pE^M_G q \iff f(p)E^N_L f(q)$ follows from Fact \ref{regular action}(i). By \ref{Borel section}(ii), we are finished.

\hfill $\blacksquare$

Note that by virtue of Proposition 2.11, Conjecture \ref{main conjecture} includes the statement that for any definable group $G$, either $G^{00} = G^{000}$ or $G^{00}/G^{000}$ is non-smooth,
which we state as a conjecture.
\begin{main conjecture}\label{special case of main conjecture}
Suppose $G$ is a group definable (or $*$-definable) over $\emptyset$ in a monster model. If $G^{000} \ne G^{00}$, then $E_0 \leq_B E_G$. 
\end{main conjecture}

\section{Computing some Borel cardinalities and Conjecture \ref{main conjecture}}

In this section (the main point of the paper), we verify versions of Conjecture \ref{main conjecture} in various situations, 
and compute some Borel cardinalities in the (very few) known examples of non $G$-compact theories. In all these examples, coverings of topological spaces and groups appear, and this seems to be more than accidental. \\

We start by showing that in the first example of a non-$G$-compact theory, described in \cite{CLPZ}, $E_L$ (on a suitable sort) has complexity $\ell^{\infty}$.  

\begin{example} We recall the many-sorted version of the example from  Section 4 of \cite{CLPZ}. We describe the standard model which has (disjoint) ``sorts" $M_{n}$ for $n = 1,2,3,\dots$. $M_{n}$ is the circle, of radius $1$, centre the origin, equipped with the clockwise oriented circular strict ordering $S_{n}(-,-,-)$ and the function $g_{n}$ which is rotation clockwise by $2\pi/n$ degrees. Let $R_{n}(x,y)$ be defined as $x=y \vee S_{n}(x,y,g_{n}(x))\vee S_{n}(y,x,g_{n}(y))$, i.e. the length of the shortest arc joining $x$ and $y$ is $< 2\pi/n$. Let $d_{n}(x,y)$ be the smallest $k$ such that there are $x=x_{0}, x_{1},\dots,x_{k} = y$ such that $R_{n}(x_{i},x_{i+1})$ for each $i=0,\dots,k-1$ (makes sense too in a saturated model). Let $M$ be the many-sorted structure $(M_{n})_{n}$ and $M^{*} = (M_{n}^{*})_{n}$ be a saturated model. We take our ``sort" $X$ to be that of infinite tuples $(a_{n})_{n}$, where $a_{n}\in M_{n}^{*}$ for each $n$. Then $X$ is a single $E_{KP}$-class, but $E_{L}((a_{n})_{n}, (b_{n})_{n})$ if and only if there is some $k<\omega$ such that
$d_{n}(a_{n},b_{n}) < k$ for all $n$.
\end{example}

\begin{proposition} In Example 3.1, $E_{L}$ above is Borel equivalent to $\ell^{\infty}$ (in the sense of the introduction). 
\end{proposition}
{\em Proof.} To do the computation, we should choose a countable model $M^{0}$ say. Let us take $M^{0}$ to be the countable elementary substructure of 
the standard model such that for each $n$, $(M^{0})_{n} $ (which we identify notationally with $(M_{n})^{0}$) consists of the  elements of $M_{n}$ 
with polar coordinates of the form $(1, 2\pi/k)$ for $k=1,2,\dots$. So the $K_{\sigma}$ equivalence relation we are considering is $E^{M^{0}}_{L}$ on 
$S_{X}(M^{0})$. Note that $S_{X}(M^{0})$ is the product of the $S_{1}((M^{0})^{n})$.  Let us fix $n$. Let $a\in (M_{n})^{*}$. Then $tp(a/(M_{n})^{0})$ 
is a ``cut" in $(M_{n})^{0}$ and hence determines a unique element of $M_{n} = {\mathbb S}^{1}$ (the circle), which is just the ``standard part" of $a$.  It is not hard to see that the equivalence relation saying that 
$a,b$ determine the same element of $\mathbb{S}^{1}$ is type-definable (over $M^{0}$) and hence that corresponding map $\pi_{n}$, say from $S_{1}((M^{0})_{n})$ 
to ${\mathbb S}^{1}$, is continuous and surjective. Hence, the map $\pi = (\pi_{n})_{n}$ from $S_{X}(M^{0})$ to the product of ${\mathbb N}$ many copies of ${\mathbb S}^{1}$ is also  continuous and surjective. It is then clear that for $p, q \in S_{X}(M^{0})$, $E^{M^{0}}_{L}(p,q)$ iff for some (any) realizations  $(a_{n})_{n}$, $(b_{n})_{n}$ of $p, q$, respectively, $E_{L}((a_{n})_{n}, (b_{n})_{n})$ iff (using the last sentence in the exposition of the example above) for some finite $k$, 
$|\pi_{n}(p) - \pi_{n}(q)| \leq 2\pi k/n$
for all $n$.  The last equivalence relation on $({\mathbb S}^{1})^{\mathbb N}$, namely the one identifying $(x_n)$ and $(y_n)$ from $({\mathbb S}^1)^{\mathbb N}$ iff for some $k$, $|x_n-y_n|\leq 2\pi k/n$ for all $n$, 
is seen to be Borel equivalent to 
$\ell^{\infty}$ as follows. It is easy to check that this equivalence relation is $K_\sigma$, so it is Borel reducible to $\ell^\infty$. On the other hand, the continuous function 
\[
{\mathbb R}^{\mathbb N}\ni (t_m)_m\mapsto \left( \exp(2\pi i \frac{t_{\rho(n)}}{n})\right)_n\in ({\mathbb S}^1)^{\mathbb N}, 
\]
where $\rho \colon {\mathbb N}\to {\mathbb N}$ is such that $\rho^{-1}(n)$ is infinite for each $n\in {\mathbb N}$, is easily checked to provide a reduction in the 
opposite direction. 

By 1.7(ii), we have proved the proposition.
\hfill $\blacksquare$\\

We now build a closely related (new) example, where the Borel complexity is $E_{0}$. 

\begin{example}  We modify the many-sorted structure $M$ = $(M_{n})_{n}$ from Example 3.1  (and from section 4 of \cite{CLPZ}).
Whenever  $n$ is an integer multiple of $m$, say $n = km$, let us add a symbol $h_{m,n}$ for the ``multiplication by $k$" map from  $M_{n}$ to $M_{m}$ (identifying both with the circle with the group operation on the circle being considered additively). Let $M'$ be $M$ equipped with these new functions. Passing to a saturated model $M'^{*}$, let $X$ now denote the set of sequences $(a_{n})_{n}$, where each $a_{n}\in M_{n}^{*}$ and where  $h_{m,n}(a_{n}) = a_{m}$ whenever $m$ divides $n$. 
\end{example}

\begin{proposition} In Example 3.3, $E_{L}$ (on $X$) is Borel equivalent to $E_{0}$. 
\end{proposition}

\noindent {\em Proof.} Let us note that for $n = km$, $h_{m,n}$ takes the function $g_{n}$ to the function $g_{m}$. 
Choose a countable model $M'^{0}$ as in the proof of Proposition 3.2, and we have again
 the map $\pi = (\pi_{n})_{n}$ from $S_{X}(M'^{0})$ to $({\mathbb S}^1)^{\mathbb N}$, but with the new choice of $X$ the image of $\pi$ is the inverse limit of the directed family of copies of ${\mathbb S}^{1}$ with respect to the various covering maps. We call this ${\widehat{\mathbb S}^{1}}$. It is a compact group. The subgroup consisting of tuples $(c_{n})_{n}$ such that $c_{n}$ converges to $0$ in ${\mathbb S}^{1}$ is a dense subgroup which we can identify with ${\mathbb R}$.

It also remains true that $X$ is a single $E_{KP}$-class and that 
that for $(a_{n})_{n}$, $(b_{n})_{n}$ in $X$, $E_{L}((a_{n})_{n}, (b_{n})_{n})$ iff for some $k$, $d_{n}(a_{n},b_{n})<k$ for all $n$. Hence,
$E_{L}((a_{n})_{n}, (b_{n})_{n})$ iff for some $k$, 
$|\pi_{n}(tp(a_{n}/M'_0)) - \pi_{n}(tp(b_{n}/M'_0))|< (2\pi k)/n$ for all $n$. 
Now, on sequences  $(c_{n})_{n}$, $(d_{n})_{n}$ from ${\widehat{\mathbb S}^{1}}$, the equivalence relation:  
for some $k$, $|c_{n}- d_{n}|< (2\pi k)/n$ for all $n$, is easily seen to be the same as being in the same coset modulo ${\mathbb R}$.  Hence, again using 
1.7, we will be done once we see that ${\widehat{{\mathbb S}^{1}}}/{\mathbb R}$ is Borel equivalent to $E_0$. That $E_0\leq_B {\widehat{{\mathbb S}^{1}}}/{\mathbb R}$ 
was first proved in \cite[Section 4]{Ro}; it also follows from \cite[Theorem 3.4.2 or Corollary 3.4.6]{BeKe}. On the other hand, since the equivalence relation 
${\widehat{{\mathbb S}^{1}}}/{\mathbb R}$ is the orbit equivalence relation of a continuous action of $\mathbb R$, there is a reduction ${\widehat{{\mathbb S}^{1}}}/{\mathbb R}\leq_B E_0$ 
by a combination of results due to Wagh and Dougherty--Jackson--Kechris; see also \cite[Theorems 1.5 and 1.16]{JaKeLo}. \hfill $\blacksquare$\\

We now pass to the examples coming from definable groups. The two basic examples from \cite{CoPi1} were:
\begin{enumerate}
\item[(i)] where $G$ is (a saturated model) of 
the universal cover of $SL_{2}({\mathbb R})$, 
\item[(ii)] where $G$ is a (saturated model of) a certain extension of $SL_{2}({\mathbb R})$ by 
$SO_{2}({\mathbb R})$ and is actually a semialgebraic group. 
\end{enumerate}
In \cite{CoPi1}, we already stated that in case (i), $G^{00}/G^{000}$ is ``naturally isomorphic to" ${\widehat{\mathbb Z}}/{\mathbb Z}$,  and in case (ii), it is ``naturally isomorphic" to 
${\mathbb S}^{1}/\Lambda$, where $\Lambda$ is a dense cyclic subgroup. 
So here we just want to check that under definitions from Section 2 these ``natural isomorphisms" can be chosen to come from or yield Borel bi-reductions (and in fact more). 

\begin{example}\label{first example from CP} In the first version of (i), the standard model is $N = (({\mathbb Z},+), ({\mathbb R},+,\cdot))$, where we feel free to add constants for finitely many elements. The relevant definable group is the universal cover $\widetilde{SL_{2}(\mathbb R)}$ of $SL_{2}({\mathbb R})$ viewed as $\{(m,g): m \in {\mathbb Z}, g\in SL_{2}({\mathbb R})\}$, where the group operation is given by 
$(m_{1},g_{1})*(m_{2},g_{2}) = (m_{1} + m_{2} + h(g_{1},g_{2}), g_{1}g_{2})$ with $h$ being a certain well-known $2$-cocycle from $SL_{2}({\mathbb R})\times SL_{2}({\mathbb R})$ to ${\mathbb Z}$  (with values $0,1, -1$ and ``definable" in the structure $N$).
We pass to a saturated model $\C$ say, of the form $(\Gamma, K)$, and let $G$ be the interpretation of the formulas defining $\widetilde{SL_{2}(\mathbb R)}$ in this big model, namely $\Gamma \times SL_{2}(K)$ with the same definition $*$ of multiplication.
\end{example}

\begin{proposition} In Example 3.5, $G = G^{00}$. Moreover, fixing a countable elementary substructure $N_{0}$ of $\C$, there is a surjective continuous
map $f$ from $S_{G}(N_{0})$ to ${\widehat{\mathbb Z}}$  (the profinite completion of ${\mathbb Z}$) which induces a bijection, in fact isomorphism of groups, between 
$S_{G}(N_{0})/E_{G}^{N_{0}}$ and ${\widehat{\mathbb Z}}/{\mathbb Z}$. Hence, by 1.7, $G^{00}/G^{000}$ is Borel equivalent to ${\widehat{\mathbb Z}}/{\mathbb Z}$  (which is, as above, known to be Borel equivalent to $E_{0}$). 
\end{proposition}
{\em Proof.}  We have the exact sequence 
$$1 \to \Gamma \to_{\iota} G \to SL_{2}(K)\to 1,$$
where the embedding $\iota: \Gamma\to G$ takes  $\gamma$ to $(\gamma, 1) \in (\Gamma\times SL_{2}(K),*) = G$.
This uses the fact that $h(1,g) = h(g,1) = 0$ for all $g\in SL(2,K)$, in particular $h(1,1) = 0$. So we will identify $\Gamma$ with its image under $\iota$. Likewise for $\mathbb{Z}$. 

Note that 
\newline
(i) $[G,G]$ (commutator subgroup of $G$) maps onto $SL_{2}(K)$. 

We also know that 
\newline
(ii) $\Gamma/\Gamma^{0} = {\widehat{\mathbb Z}}$. 
\newline
This equality is as topological groups, where the topology on the first is the logic topology. Note that the map $\Gamma \to \Gamma/\Gamma^{0}$ factors through $S_{\Gamma}(N_{0})$ and that this logic topology coincides with the quotient topology with respect to $S_{\Gamma}(N_{0})$ 

Moreover, it  was proved in \cite{CoPi1} (proof of Theorem 3.2 there) that:  
\newline
(iii) $[G,G]\cap \Gamma = \mathbb{Z}$ (which note is disjoint from $\Gamma^{0}$),
\newline
(iv) $G^{000} = \Gamma^{0}\cdot [G,G]$,
\newline
(v) $G^{000}=(\Gamma^0 + \mathbb{Z}) \times SL_2(K)$ as sets,
\newline
(vi) $G^{00} = G$. 

It follows from the above that $\iota$ induces an isomorphism of groups which we still call $\iota: \Gamma/(\Gamma^{0} + {\mathbb Z}) \to G/G^{000}$. Note that via (ii) above, $\Gamma/(\Gamma^{0} + {\mathbb Z})$ can be identified with  
${\widehat {\mathbb Z}}/{\mathbb Z}$.   
It is not so hard to see that $\iota$ is induced by a Borel function from ${\widehat{\mathbb Z}}$ to
$S_{G}(N_{0})$, yielding  ${\widehat{\mathbb Z}}/\mathbb{Z}  \leq_{B} S_{G}^{N_{0}}/E_{G}^{N_{0}}$.  But we need more. There are various options, and here is one of them.   

We will identify $\Gamma/\Gamma^{0}$ and ${\widehat{\mathbb Z}}$ via (ii) above and the remarks following it. 
Let $f: S_{G}(N_{0}) \to {\widehat{\mathbb Z}}$  be given by 
$f(tp((\gamma,g)/N_{0})) = \gamma/\Gamma_{0}$.  The map $f$ is well-defined and clearly continuous and surjective.  We will identify, via Remark 2.9(iii), $S_{G}(N_{0})/E_{G}^{N_{0}}$ and $G/G^{000} = G^{00}/G^{000}$.\\[3mm] 
{\bf Claim} $(\gamma,g)$ and $(\gamma', g')$ are in the same coset modulo $G^{000}$ if and only if $\gamma/\Gamma^{0}$ and $\gamma'/\Gamma^{0}$ are in the same coset modulo $\Gamma^{0}+ {\mathbb Z}$. Moreover, the map taking $(\gamma,g)/G^{000}$ to $\gamma/(\Gamma^{0}+{\mathbb Z})$ yields an isomorphism of groups.\\[3mm] 
%
%
{\em Proof of claim.} By the description of $G^{000}$ in (v), the formula for multiplication in $G$ and the fact that the 2-cocycle $h$ takes values in $\mathbb{Z}$, we see that if $(\gamma,g)$ and $(\gamma',g')$ are in the same coset of $G^{000}$, then $\gamma = \gamma' + \gamma_{1} + a$ for some $a\in {\mathbb Z}$ and $\gamma_{1}\in \Gamma^{0}$. The same argument proves the converse, yielding  the first part of the claim. This together with the formula for multiplication in $G$ shows that $(\gamma,g)/G^{000} \mapsto \gamma/(\Gamma^{0}+{\mathbb Z})$ yields a well-defined isomorphisms between  $G^{00}/G^{000}$ and $\Gamma/(\Gamma^0 + \mathbb{Z})$.\\

The Claim finishes the proof of Proposition 3.6. \hfill $\blacksquare$\\

Note that as $(SL_{2}({\mathbb R}), \cdot)$ is bi-interpretable with $({\mathbb R}, +, \cdot)$, the structure $N$ in Example 3.5 is essentially the $2$-sorted structure
$(({\mathbb Z}, +), (SL_{2}({\mathbb R}), \cdot))$, which as we saw above interprets 
naturally the group $\widetilde{SL_{2}({\mathbb R})}$. 
However, 
$(({\mathbb Z}, +), (SL_{2}({\mathbb R}), \cdot))$ and $(\widetilde{SL_{2}({\mathbb R})}, \cdot)$ are not bi-interpretable even with parameters. So we ask what goes on in this latter reduct (where information is lost). 
Let $M$ be the structure  
$(\widetilde{SL_{2}({\mathbb R})}, \cdot)$. 

\begin{proposition} Let $(G,\cdot)$ be a saturated model of $Th(M)$, and $M_{0}$ a countable model. Then $G^{00} = G$, and $G/G^{000}$ is again Borel  equivalent to 
$\widehat{\mathbb Z}/{\mathbb Z}$.  
\end{proposition}
{\em Proof.} We may assume that $(G,\cdot)$ is the group from Example 3.5,  definable in the monster model $\C$ from there.  (So $(G,\cdot)$ is a reduct of $\C$.) Let $N_{0}$ be a countable elementary substructure of $\C$ as in 3.6 and let $M_{0}$ be its reduct (in the language of groups). 

Note that the natural surjective homomorphism $\pi: G \to SL_{2}(K)$  ($K$ a saturated real closed field) is definable. We let $\Gamma$ denote the kernel which we again write additively. It contains a canonical copy of ${\mathbb Z}$
(the kernel of the restriction of $\pi$ to $\widetilde{SL_{2}({\mathbb R})}$), disjoint from $\Gamma^{0}$. 
And the properties  $G = G^{00}$, $[G,G]\cap \Gamma = \mathbb{Z}$ and $G^{000} = \Gamma^{0}\cdot [G,G]$ are inherited from Example 3.5.
Fix a countable model $M_{0}$. We have the continuous surjective map from 
$S_{\Gamma}(M_{0})$ to $\Gamma/\Gamma^{0}$ (= $\widehat{\mathbb Z}$), taking 
$tp(\gamma/M_{0})$ to $\gamma/\Gamma^{0}$. By 1.7, this has a Borel section $f:\widehat{\mathbb Z} \to S_{\Gamma}(M_{0})\subset S_{G}(M_{0})$ say, and the methods from the proof of Proposition 3.6 give that $f$ is a Borel reduction from 
$\widehat{\mathbb Z}/{\mathbb Z}$ to  $S_{G}(M_{0})/E_{G}^{M_{0}}$. 

On the other hand, we have a surjective continuous function from the type space 
$S_{G}(N_{0})$  (Example 3.5 and  Proposition 3.6) to the type space 
$S_{G}(M_{0})$. By 1.7, this has a Borel section $h$ say, and composing it with the continuous surjective map $f:S_{G}(N_{0}) \to \widehat{\mathbb Z}$ from Proposition 3.6 gives the required Borel reduction from $S_{G}(M_{0})/E_{G}^{M_{0}}$ to $\widehat{\mathbb Z}/{\mathbb Z}$. This completes the proof.
\hfill $\blacksquare$

\begin{example}\label{second example from CP} Let $h: SL_{2}(\mathbb{R}) \to \mathbb{Z}$ be the $2$-cocycle from  Example 3.5.  Let us denote by $\mathbb{R}/\mathbb Z$ the semialgebraic group $[0,1)$ with addition modulo $1$. Fix an element $c \in \mathbb{R}/\mathbb{Z}$ of infinite order.  Consider the connected real semialgebraic group whose universe is
$\mathbb{R}/\mathbb{Z} \times SL_{2}(\mathbb{R})$ equipped with the operation $*$, where $(a_{1},g_{1})*(a_{2},g_{2}) = (a + b + h(g_{1},g_{2})c, g_{1}g_{2})$.  Pass to a saturated real closed field $K$, and let $G$ be the group definable in $K$ by the same formulas. 
\end{example}

\begin{proposition} $G = G^{00}$ and $G/G^{000}$ is Borel equivalent to $\mathbb{S}^{1}/\Lambda$  (in fact, to  $(\mathbb{R}/\mathbb{Z})/\langle c \rangle$), 
where $\Lambda$ is a cyclic dense subgroup. And this is again known to be Borel equivalent to $E_{0}$.  
\end{proposition}
{\em Proof.} The universe of $G$ is $\Gamma\times SL_{2}(K)$ with the group operation $*$ defined as above. 
Let the cyclic group generated by $c$ be denoted by $\Lambda$.  From Theorem 3.3 of \cite{CoPi1}, $G = G^{00}$, $G^{000} = \Gamma^{00}\cdot[G,G]$, 
and moreover, $[G,G]\cap \Gamma = \Lambda$. And $\Gamma/\Gamma^{00}$ identifies with $\mathbb{R}/\mathbb{Z}$. 

Fix a countable model $M_{0}$.
As in the proof of 3.6, the continuous surjective map from $S_{G}(M_{0})$ to $\Gamma/\Gamma^{00}$ which takes $tp((\gamma,g)/M_{0})$ to $\gamma/\Gamma^{00} \in \mathbb{R}/\mathbb{Z}$ induces a group isomorphism between $G/G^{000}$ and $(\mathbb{R}/\mathbb{Z})/\Lambda$, and by 1.7, we are finished. 
\hfill $\blacksquare$\\

In \cite{CoPi2}, $G^{00}/G^{000}$ is studied for arbitrary definable groups in (saturated) $o$-minimal expansions of real closed fields, and it is shown that either $G^{00} = G^{000}$ or $G^{00}/G^{000}$ is (abstractly) isomorphic to the quotient $A/\Gamma_{0}$ of a compact connected commutative real Lie group $A$ by a finitely generated dense subgroup $\Gamma_{0}$. A further analysis, which we will not carry out here, yields that in this latter case $G^{00}/G^{000}$ is in fact Borel equivalent to $A/\Gamma_{0}$  (and thus again to $E_{0}$). \\

We now want to use products  to  produce $*$-definable groups  (even in real closed fields), where $G^{00}/G^{000}$ has Borel complexity $\ell^{\infty}$.  We first use the {\em finite covers} of $SL_{2}(\mathbb{R})$ to answer a question of Gismatullin \cite{Gi}.  Let us fix a group $G$ defined or even $*$-defined over $\emptyset$ say. 
Recall from Fact 1.10 that $G^{000}$ (which by convention means $G^{000}_{\emptyset}$) is the subgroup of $G$ generated by $\{ab^{-1}: a,b\in G$ and $E_{L}(a,b)\}$. As $E_{L}$ is the transtive closure of the partial type $\Theta(x,y)$ (expressing that $x,y$ begin an infinite indiscernible sequence), it follows that $G^{000}$ is generated by $X_{\Theta} = \{ab^{-1}: a,b\in G$ and $\Theta(a,b)\}$. If $G^{000}$ is generated in finitely many steps by $X_{\Theta}$, then clearly $G^{000}$ is type-definable and equals $G^{00}$, and in this case we call the minimum number of steps required the ``diameter" of $G^{00}$. 
Conversely, if $G^{000}$ is type-definable, then $G^{000}=G^{00}$ has finite diameter \cite[Remark 3.8]{GiNe} (but this is much more involving as it uses \cite{Ne}). The question of Gismatullin was whether for any $n$ there is a group $G_{n}$ such that $G_{n} = G_{n}^{000}$ and $G$ has (finite) diameter $\geq n$.   (The idea being that as in Example  3.1 this would enable one to produce $*$-definable groups with $G^{00}\neq G^{000}$.) Note that in theories where $G^{000}_{A}$ does not depend on the choice of $A$, we can  work over any set of parameters, such as a model $M$. But then types over $M$, KP types over $M$ and Lascar strong types over $M$ coincide. And so $G^{000}$ is generated by 
$\{ ab^{-1}: tp(a/M) = tp(b/M) \}$  (which is type-definable and roughly corresponds to $\Theta(x,y)$ over $M$), and we may as well consider  diameter from this point of view. 

In any case, this lack of a bound on finite diameters is precisely what is behind Example 3.5.  The various finite covers of $SL_{2}(\mathbb{R})$ (looked at in a saturated model) are groups $G$ such that $G = G^{000}$ but have arbitrarily large finite diameters, as we explain now.  
Let $h:SL_{2}(\mathbb{R})\times SL_{2}(\mathbb{R}) \to \mathbb{Z}$ be the cocycle from 3.5 which induces the group operation $*$ on $\mathbb{Z}\times SL_{2}(\mathbb{R})$ giving $\widetilde{SL_{2}(\mathbb{R})}$. As before, we identify the group $(\mathbb{Z},+)$ with $\{(a,1):a\in \mathbb{Z}\}$. 
Then 
for any $n$ the quotient of $\widetilde{SL_{2}(\mathbb{R})}$ by the central subgroup $n\mathbb{Z}$ identifies with
$(\mathbb{Z}/n\mathbb{Z}) \times SL_{2}(\mathbb{R})$ with group operation given by the same formula, namely
$(k_{1}, g_{1})*(k_{2}, g_{2}) = (k_{1} + k_{2} + h(g_{1}, g_{2}), g_{1}g_{2})$, where 
$k_{1} + k_{2} + h(g_{1}, g_{2})$ is computed in $\mathbb{Z}/n\mathbb{Z}$. (Remember that $h$ has values in $\{-1,0,1\}$.) So this finite cover of 
$SL_{2}(\mathbb{R})$  is a definable group in $(\mathbb{R},+,\times)$. Let $G_{n}$ be its interpretation in a saturated real closed field 
$(K,+,\times)$. So the $G_{n}$ are (definably connected) finite central extensions of $SL_{2}(K)$.

With this notation:
\begin{lemma} $G_{n} =  G_{n}^{000}$, but has finite diameter $\geq O(n)$. 
\end{lemma}
{\em Proof.} Note that $G_{n}$ is connected (no proper definable subgroup of finite index). As $G_{n}$ projects definably 
onto $SL_{2}(K)$ with finite kernel $\mathbb{Z}/n\mathbb{Z}$, and $SL_{2}(K)$ has no proper non-trivial normal subgroups, we see that $G_{n} = 
G_{n}^{000}$. 

Now, we consider the diameters. We are working in the theory $RCF$ where $dcl(\emptyset)$ is an elementary substructure. So, if
$a =(k_{1}, g_{1})$ and $ b = (k_{2},g_{2})$ are in $G_{n}$ such that either $\Theta(a,b)$ or even $a$ and $b$ have the same type, then $k_{1} = k_{2}$. 
So, $b^{-1} = (-k_{1} - h(g_{2}^{-1},g_{2}), g_{2}^{-1})$. Thus, $a*b^{-1} = 
(-h(g_{2}^{-1},g_{2}) + h(g_{1},g_{2}^{-1}), g_{1}g_{2}^{-1})$, and the first coordinate is between $-2$ and $2$. 
Hence, any product of $r$ such 
elements has first coordinate in between $-3r$ and $3r$ (computed mod $n$). So the product of at least $O(n)$ such elements is needed to cover 
$G_{n}^{*}$.  We also need to know that $G_{n}$ has finite diameter. 
It follows from the fact that $G_n^{000}=G_n$ and \cite[Remark 3.8]{GiNe}, but it can also be seen more directly. Namely, it is well-known that $SL_{2}(K)$ has finite diameter. Also, as $G_{n} = G_{n}^{000}$, some finitely many (finite) products of elements $ab^{-1}$, where $\Theta(a,b)$ holds,  will cover the (finite) kernel $\mathbb{Z}/n\mathbb{Z}$. Hence, $G_{n}$ has finite diameter.
\hfill $\blacksquare$\\

We will now let $\Gamma$ denote $\prod_{n}\mathbb{Z}/n\mathbb{Z}$. $\Gamma$ is a so-called $*$-definable group in $(K,+,\times)$. As it is also ``bounded", it is a compact topological group when equipped with the logic topology, and this topology of course coincides with the natural product topology on $\Gamma$. It is also of course separable. We call a sequence $(a_{n})_{n}\in \prod_{n}\mathbb{Z}/n\mathbb{Z}$ {\em bounded} if for some natural number $r$, for all $n$, $a_{n}$ is in (or has a representative in) the interval $(-r,r)$. The collection of such bounded sequences is clearly a countable dense subgroup of $\Gamma$ which we call $B_{\Gamma}$.  The quotient $\Gamma/B_{\Gamma}$ is known to be Borel equivalent to $\ell^{\infty}$.

\begin{proposition} Let $G$ be the product $\prod_{n}G_{n}$. So $G$ is a so-called $*$-definable group in $(K,+,\times)$. Then $G^{00} = G$, $G^{000}$ is the commutator subgroup $[G,G]$, and $G/G^{000}$ is Borel equivalent (and isomorphic as a group) to $\Gamma/B_{\Gamma}$ (from the above paragraph). 
\end{proposition}
{\em Proof.}  Let $R$ denote $\prod_{n}SL_{2}(K)$. We have the surjective ($*$-definable) homomorphism $\pi:G\to R$ with kernel $\Gamma$ (central in $G$). We go through several claims.\\[3mm]
{\bf Claim 1.} $R = R^{000}$.\\[3mm] 
%
{\em Proof of Claim 1.}  
This follows from \cite[Theorem 3.11, Lemma 2.12(1), Theorem 2.15(2)]{Gi} (as was told us by J. Gismatullin).
\newline
But here is another direct proof. First find some $g = (g_{n})_{n}\in R^{000}$ such that no $g_{n}$ is central in 
$SL_{2}(K)$. (If there is no such $g$, then let $h_{\alpha}$ for $\alpha < \kappa$ be elements of $SL_{2}(K)$ in 
different cosets modulo the centre. For each $\alpha < \kappa$, let  $g^{\alpha}\in R$ have $nth$ coordinate $h_{\alpha}$ for each $n$. 
Then the $g^{\alpha}$ are in different cosets modulo $R^{000}$. As $\kappa$ is arbitrary this contradicts bounded 
index of $R^{000}$). Now clearly the conjugacy class of $g$ in $R$ generates $R$ (as it is true in each coordinate).
Hence as $R^{000}$ is normal in $R$, $R = R^{000}$. 
\hfill $\blacksquare$\\

As $SL_{2}(K)$ is perfect (of finite commutator width), we get\\[3mm]
{\bf Claim 2.} $R$ is perfect (of finite commutator width).\\[3mm]
As $G$ is a central extension of a perfect group, we obtain\\[3mm]
{\bf Claim 3.} $[G,G]$ is perfect, and moreover maps onto $R$ under $\pi$.\\[3mm]
%
Hence, as $[G,G]$ is invariant and  has bounded index in $G$, we get\\[3mm]
{\bf Claim 4.} $G^{000} \leq [G,G]$.\\[3mm]
%
{\bf Claim 5.}  $G=\Gamma \cdot G^{000}$ and $G^{000} = [G,G]$.\\[3mm]
%

{\em Proof of Claim 5.} By Claim 1, $\pi[G^{000}]=R^{000}=R$, hence $G=\Gamma \cdot G^{000}$. Moreover, $\Gamma=ker(\pi)$ is central in $G$. Thus, $[G,G]=[\Gamma \cdot G^{000},\Gamma \cdot G^{000}]=[G^{000},G^{000}]\leq G^{000}$. So, the second part of Claim 5 follows from Claim 4. \hfill $\blacksquare$\\

As $G = \Gamma\cdot G^{000}$, it follows that\\[3mm] 
{\bf Claim 6.}  $\Gamma/(\Gamma\cap G^{000})$ is isomorphic to $G/G^{000}$ (induced by the embedding of $\Gamma$ in $G$).\\

It remains to identify $\Gamma\cap G^{000}$ as the subgroup $B_{\Gamma}$ and to show that we
have a Borel equivalence in addition to an isomorphism. 

For the first part, we look at our representation of $G_{n}$ as 
$(\mathbb{Z}/n\mathbb{Z} \times SL_{2}(K),*)$.  The natural embedding taking $(a_{n})_{n}\in \Gamma$ to $(a_{n},1)_{n}\in G$ is a group embedding. 
As $G^{000} = [G,G]$, clearly $(a_{n})_{n}$ is in $\Gamma\cap G^{000}$ iff for some $r$, $(a_{n})_{n}$ is a product  of at most $r$ commutators in $G$ iff for for some $r$ and  each $n$, $a_{n}$ is a product of at most $r$ commutators in $G_{n}$. 
Using the formula for multiplication in $G_n$ and the fact that the 2-cocycle $h$ has finite image, this shows that $\Gamma \cap G^{000} \leq B_\Gamma$.  Using additionally the fact that $\mathbb{Z}\leq [\widetilde{SL_2(R)}, \widetilde{SL_2(R)}]$ (which follows from the perfectness of $\widetilde{SL_2(R)}$), this also implies the opposite inclusion, i.e. $B_\Gamma \leq \Gamma \cap G^{000}$. So, we have proved the following\\[3mm]
%
{\bf Claim 7.} $\Gamma\cap G^{000} = B_\Gamma$.\\ 

At this point one sees that $G = G^{00}$. For $G^{00}$ contains $G^{000}$, hence $G^{00}\cap \Gamma$, a closed subgroup of $\Gamma$, contains, by Claim 7, the dense subgroup $B_\Gamma$, so $G^{00}$ contains $\Gamma$, 
and so equals $G$ by Claim 5.\\[3mm]
%
%
{\bf Claim 8.} $G^{000} = \{(a_{n},g_{n})_{n}: (a_{n})_{n} \in B_\Gamma, g_{n}\in SL_{2}(K) \; \mbox{for all} \;n \}$, which can be identified with $B_\Gamma \times R$.\\[3mm]
{\em Proof of Claim 8.}  Since  the 2-cocycle  $h$ has finite image, we easily get that the set $ B_\Gamma \times R$ forms a subgroup of $G$. This subgroup is of bounded index and it is invariant by Claim 7. Therefore,  $G^{000} \leq B_\Gamma \times R$. As $\pi[G^{000}] =R$ (by Claim 1) and $\Gamma\cap G^{000} = B_\Gamma$ (by Claim 7), we conclude that $G^{000} = B_\Gamma \times R$.\hfill $\blacksquare$\\

Let  $M_{0}$ be $dcl(\emptyset)$, a countable elementary substructure of $(K,+,\times)$. In particular, each $\mathbb{Z}/n\mathbb{Z}$ is contained in $M_{0}$. 
Let $f:S_{G}(M_{0}) \to \Gamma$  take $tp((a_{n},g_{n})_{n}/M_{0})$ to $tp((a_{n})_{n}/M_{0})$ (= $(a_{n})_{n} \in \Gamma$ ). It is clearly a continuous surjection from $S_{G}(M_{0})$ to $\Gamma$. It follows from Claim 8 (and the fact that the $2$-cocycle $h$ had image $\{-1,0,1\}$) that  $(a_{n},g_{n})_{n}*(b_{n},h_{n})_{n}^{-1} \in G^{000}$ if and only if
$(a_{n}-b_{n})_{n} \in B_{\Gamma}$.  By Fact 1.7 (ii), $f$ yields the Borel 
equivalence of $G/G^{000}$ and $\Gamma/B_{\Gamma}$. In fact, note that the bijection 
between $G/G^{000}$ and $\Gamma/B_{\Gamma}$ induced by $f$ is the natural isomorphism 
of groups mentioned in Claim 6.

\hfill $\blacksquare$\\

\begin{remark} Instead of taking the product of the $G_{n}$'s we could take the inverse limit under the natural surjective homomorphisms $G_{n} \to G_{m}$ when $m$ divides $n$.  Let us call this group $H$. Then $H$ projects onto $SL_{2}(K)$ with kernel $\widehat{Z}$. And one can prove in a similar fashion to the above that
$H = H^{00}$ and $H/H^{000}$ is Borel equivalent to $\widehat{\mathbb{Z}}/\mathbb{Z}$ (which we know has complexity $E_{0}$).

\end{remark}



 \section{$\pmb{G^{000}}$ and Borel reductions in a general context}

In \cite{GiKr}, a certain axiomatic framework for understanding the group examples discussed in the previous section was given, leading to some new classes of examples of definable groups $G$ with $G^{000}\neq G^{00}$ and also raising several questions. More precisely, \cite{GiKr} deals with the question when an extension $\widetilde{G}$ of a group $G$ by an abelian group $A$ given by a 2-cocycle $h$ with finite image satisfies $\widetilde{G}^{000} \ne \widetilde{G}^{00}$ (in a monster model). Paper \cite{GiKr} provides a sufficient condition on $h$ to have $\widetilde{G}^{000} \ne \widetilde{G}^{00}$, which also turns out to be necessary in a rather general situation. 

In this section, we study Conjecture \ref{main conjecture} (more precisely, its definable group version -- Conjecture \ref{special case of main conjecture}) in the context from \cite{GiKr}, giving several positive results (see Corollary \ref{conjecture 1 true}), but also giving a possible scenario which leads to a counterexample (see Proposition 4.7). Corollary \ref{conjecture 1 true} implies that Conjecture \ref{special case of main conjecture} holds for classes of concrete examples obtained in \cite[Section 4]{GiKr}, in particular, for Examples \ref{first example from CP} and \ref{second example from CP} described in Section 3 (which was proved in Section 3 more directly).

This section is somewhat technical, and some familiarity with the preprint \cite{GiKr} would help the reader. In this section, we write the parameter sets explicitly (even if it is $\emptyset$).
Before recalling the situation from \cite{GiKr}, let us make a certain observation, which will be very useful later.

\begin{proposition}\label{E_0 < E_G}
Let $G$ be a group $\emptyset$-type-definable in a monster model $\C$ (even on infinite tuples).
Suppose that for some countable set $B$ there exits a $B$-type-definable subgroup $H$ of $G$ such that $H^{00}_B \leq G^{000}_\emptyset \cap H$. Then $E_a \leq_B E_G$, where $E_a$ is the relation on $(G^{00}_\emptyset \cap H)/H^{00}_B$ (equipped with the logic topology) of lying in the same coset modulo $(G^{000}_\emptyset \cap H)/H^{00}_B$. If moreover $G^{000}_\emptyset \cap H$ is not type-definable, then  $E_0 \leq_B E_a$, and so $E_0 \leq_B E_G$.
\end{proposition}
{\em Proof.} Recall that $H/H^{00}_B$ can be equipped with the so-called logic topology (whose closed sets are the sets with type-definable preimages under the quotient map $H \to H/H^{00}_B$). In this way, $H/H^{00}_B$  becomes a compact topological group. Since both the language and the set $B$ are countable, this is a Polish group. 

Put $H_1=G^{00}_\emptyset \cap H$. Then, $H^{00}_B \leq G^{000}_\emptyset \cap H \leq H_1 \leq H$, and $H_1/H^{00}_B$ is a closed subgroup of $H/H^{00}_B$, so it is also a compact, Polish group.

Choose a countable $M\prec \C$ containing $B$. Recall that $S_{H_1}(M):=\{ tp(h/M): h \in H_1\}$ is a compact, Polish space . Define a binary relation $E$ on $S_{H_1}(M)$ by 
$$pEq \iff (\exists a \models p)(\exists b \models q) (ab^{-1} \in H^{00}_B).$$
It is easy to see that that $pEq \iff (\forall a \models p)(\forall b \models q) (ab^{-1} \in H^{00}_B)$ (see the proof of Remark \ref{equivalence relation}), and so $E$ is a closed equivalence relation on $S_{H_1}(M)$. Thus, the quotient space $S_{H_1}(M)/E$ is a compact, Polish space.

Let $f \colon H_1/H^{00}_B \to S_{H_1}(M)/E$ be defined by $f(hH^{00}_B)=[tp(h/M)]_E$. It is clear that $f$ is a homeomorphism.

Let $E_a$ be the equivalence relation on $H_1/H^{00}_B$ of lying in the same coset modulo $(G^{000}_\emptyset \cap H)/H^{00}_B$. By Fact \ref{G^{000}} and the definition of the logic topology, we get that $E_a$ is of type $K_\sigma$ and so Borel.\\[3mm]
{\bf Claim 1.} $E_a \leq_B E_G$.\\[3mm]
{\em Proof of Claim 1.} By Fact \ref{Borel section}, the quotient map from $S_{H_1}(M)$ onto $S_{H_1}(M)/E$ has a Borel section $f_1\colon S_{H_1}(M)/E \to S_{H_1}(M)$. We will show that the composition $f_1 \circ f \colon H_1/H^{00}_B \to S_{H_1}(M) \subseteq S_{G^{00}_\emptyset}(M)$ is a Borel reduction of $E_a$ to $E_G$.

Since $f$ and $f_1$ are Borel, so is $f_1 \circ f$. It remains to show that for any $h_1,h_2 \in H_1$
$$h_1h_2^{-1} \in G^{000}_\emptyset \iff f_1([tp(h_1/M)]_E)\, E_G \, f_1([tp(h_2/M)]_E).$$

We have that 
$$h_1h_2^{-1} \in G^{000}_\emptyset \iff tp(h_1/M)\, E_G \, tp(h_2/M),$$
so we will be done if we notice that
$$tp(h_1/M)\, E_G \, f_1([tp(h_1/M)]_E) \;\, \mbox{and}\;\, tp(h_2/M)\, E_G \, f_1([tp(h_2/M)]_E).$$
For this, consider any $p,q \in S_{H_1}(M)$ such that $p E q$. By the definition of $E$, there are $a \models p$ and $b \models q$ with $ab^{-1}\in  H^{00}_B$. Then $ab^{-1} \in G^{000}_\emptyset$, so $pE_G q$. \hfill $\blacksquare$\\[3mm]
{\bf Claim 2.} $E_0 \leq_B E_a$ (assuming that $G^{000}_\emptyset \cap H$ is not type-definable).\\[3mm]
{\em Proof of Claim 2.}
Let $H_0$ be the preimage of the closure of $(G^{000}_\emptyset \cap H)/H^{00}_B$ by the quotient map from $H_1$ onto $H_1/H^{00}_B$. Then $(G^{000}_\emptyset \cap H)/H^{00}_B$ acts on $H_0/H^{00}_B$ by translations, which are of course homeomorphisms. Moreover, the orbit of the neutral element $H^{00}_B$ under this action is dense, and the relation $E_a$ restricted to $H_0/H^{00}_B$, which we denote by $E_a'$,  is still of type $K_\sigma$. So, the conclusion follows from Fact \ref{dense orbit}, provided that we prove that $E_a'$ is a meager subset of $H_0/H^{00}_B \times H_0/H^{00}_B$.

Suppose this is not the case. Since $E_a'$ is of type $K_\sigma$, this means that $E_a'$ has a non-empty interior, which easily implies that $(G^{000}_\emptyset \cap H)/H^{00}_B$ is an open subgroup of  $H_0/H^{00}_B$. So, it is also a closed subgroup, hence $G^{000}_\emptyset \cap H$ is type-definable, a contradiction. \hfill $\blacksquare$\\

The main result of \cite{GiKr}, which we recall below, is in the language of group extensions and $2$-cocycles. The reader is referred to the initial part of \cite[Section 2]{GiKr} for a very short overview of some basic notions and facts concerning these issues.
   
We consider a situation when groups $G$ and $A$ together an action $\cdot$ of $G$ on $A$ are $\emptyset$-definable in a (many-sorted) structure $\G$ (e.g. $\G$ consists of the pure groups $G$ and $(A,+)$ together with the action of $G$ on $A$).
Moreover, it is assumed in \cite{GiKr} that $h \colon G \times G \to A$ is a $B$-definable 2-cocycle with finite image contained in the finite set $B$ of parameters from $\G$. In fact, instead of $Im(h) \subseteq B$ it is enough to require that $Im(h) \subseteq dcl(B)$, so adding $B$ to the language, we can assume that $B=\emptyset$ and require $Im(h) \subseteq dcl(\emptyset)$.

Let $\widetilde{G}$ be the extension of $G$ by $A$ defined by means of $h$, i.e. $\widetilde{G}$ is the product $A \times G$ with the following group law 
$$(a_1,g_1)*(a_2,g_2) = (a_1+g_1\cdot a_2+h(g_1,g_2),g_1g_2).$$

The group $\widetilde{G}$ is, of course, definable in $\G$.
Let $\G^* \succ \G$ be a monster model. Denote by $G^*$ the interpretation of $G$ in $\G^*$, by $A^*$ the interpretation of $A$ in $\G^*$, and by $\widetilde{G^*}$ the interpretation of $\widetilde{G}$ in $\G^*$. We have the following exact sequence
%
%
$$
\xymatrix{
 1 \ar@{^{(}->}[r] & A^* \ar@{^{(}->}[r] & {}\widetilde{G^*} \ar@{>>}[r]^-{\pi} & G^* \ar@{>>}[r] & 1},
$$
where $\pi$ is the projection on the second coordinate.

In the sequel, by ${A^*}^{000}_\emptyset$, ${G^*}^{000}_\emptyset$ and $\widetilde{G^*}^{000}_\emptyset$ we denote the smallest subgroups of bounded index of $A^*$, $G^*$ and $\widetilde{G^*}$, respectively, which are invariant under $Aut(\G^*)$.  The components ${A^*}^{00}_\emptyset$, ${G^*}^{00}_\emptyset$ and $\widetilde{G^*}^{00}_\emptyset$ are also computed working in $\G^*$.
By $A_0$ we will denote the subgroup of $A$ generated by the image of $h$. Notice that $A_0$ is finitely generated and contained in $dcl(\emptyset)$.

The following is \cite[Theorem 2.1]{GiKr} (after adding $B$ to the language).

\begin{fact} \label{main result from GiKr}
Let $G$ be a group acting by automorphisms on an abelian group $A$, everything $\emptyset$-definable in a structure $\G$, and let $h\colon G \times G \to A$ be a 2-cocycle  which is $\emptyset$-definable in $\G$ and with finite image $Im(h)$ contained in $dcl(\emptyset)$. Let $A^*_1$ be a bounded index subgroup of $A^*$ which is type-definable over $\emptyset$ and which is invariant under the action of $G^*$. Assume that: 
\begin{enumerate}
\item[(i)] 
the induced 2-cocycle $\overline{h}\colon {G^*}^{00}_\emptyset \times {G^*}^{00}_\emptyset \to A_0/\left({A^*}_{1} \cap A_0\right)$ is non-splitting,
\item[(ii)] $A_0/\left( {A^*}_{1} \cap A_0 \right)$ is torsion free (and so isomorphic to ${\mathbb Z}^n$ for some natural $n$).  
\end{enumerate}
Then $\widetilde{G^*}^{000}_\emptyset \ne \widetilde{G^*}^{00}_\emptyset$.

Suppose furthermore that ${G^*}^{000}_\emptyset = G^*$, and for every proper, $\emptyset$-type-definable in $\G^*$ subgroup $H$ of $A^*$ with bounded index the induced 2-cocycle $\overline{h'}\colon G^* \times G^* \to A_0/\left(H \cap A_0\right)$ is non-splitting. Then $\widetilde{G^*}^{00}_\emptyset = \widetilde{G^*}$.
\end{fact}

Using this theorem, Section 4 of \cite{GiKr} provides  some new classes of examples of groups of the form $\widetilde{G^*}$ for which $\widetilde{G^*}^{000}_\emptyset \ne \widetilde{G^*}^{00}_\emptyset$ (generalizing   \cite{CoPi1}). 

Our next goal is to understand the Borel cardinality of $E_{\widetilde{G^*}}$. We will usually assume that ${G^*}^{000}_\emptyset ={G^*}^{00}_\emptyset$ and ${A^*}^{000}_\emptyset ={A^*}^{00}_\emptyset$, which means that the only reason for $\widetilde{G^*}^{000}_\emptyset \ne \widetilde{G^*}^{00}_\emptyset$ comes from the interaction between $G^*$ and $A^*$ via the action of $G^*$ on $A^*$ and the 2-cocycle $h$ (and not from the fact that already $G^*$ or $A^*$ is an example where the two components differ; notice that ${G^*}^{000}_\emptyset \ne {G^*}^{00}_\emptyset$ implies $\widetilde{G^*}^{000}_\emptyset \ne \widetilde{G^*}^{00}_\emptyset$). For example, the assumption ${A^*}^{000}_\emptyset ={A^*}^{00}_\emptyset$ holds when $\G=((A,+),M)$, where $G$ is $\emptyset$-definable in the structure $M$ and it acts trivially on $A$, because then ${A^*}^{000}_\emptyset ={A^*}^{00}_\emptyset={A^*}^0$ is the intersection of all definable in $(A^*,+)$ subgroups of $A^*$ of finite index. This situation takes place in concrete applications of Fact \ref{main result from GiKr} in Section 4 of \cite{GiKr}. The assumption ${G^*}^{000}_\emptyset ={G^*}^{00}_\emptyset$ is also satisfied in these applications.

\begin{proposition}\label{x}
Consider the situation from the first sentence of Fact \ref{main result from GiKr}. Assume that $\widetilde{G^*}^{000}_\emptyset \cap A^*$ is not type-definable and ${A^*}^{00}_\emptyset \leq \widetilde{G^*}^{000}_\emptyset \cap A^*$. Then $E_0 \leq_B E_{\widetilde{G^*}}$. 
\end{proposition}
Before the proof, notice that the assumption ${A^*}^{00}_\emptyset \leq \widetilde{G^*}^{000}_\emptyset \cap A^*$ is satisfied whenever ${A^*}^{000}_\emptyset = {A^*}^{00}_\emptyset$.\\[3mm]
{\em Proof.} Since ${A^*}^{00}_\emptyset \leq \widetilde{G^*}^{000}_\emptyset \cap A^*$ and $\widetilde{G^*}^{000}_\emptyset \cap A^*$ is not type-definable, the desired conclusion follows from Proposition \ref{E_0 < E_G}, putting $H=A^*$. \hfill $\blacksquare$

\begin{proposition}\label{xx}
Consider the situation from the first sentence of Fact \ref{main result from GiKr}. Assume that ${G^*}^{000}_\emptyset ={G^*}^{00}_\emptyset$ and there exists a bounded index subgroup $A^*_1$ of $A^*$ which is $\emptyset$-type-definable (in $\G^*$) and such that $A^*_1 \subseteq \widetilde{G^*}^{000}_\emptyset \cap A^*$. Let $F$ be the relation on the compact, Polish group $(\widetilde{G^*}^{00}_\emptyset \cap A^*)/A^*_1$ of lying in the same coset modulo $(\widetilde{G^*}^{000}_\emptyset \cap A^*)/A^*_1$. Then $E_{\widetilde{G^*}} \sim_B F$. 
\end{proposition}

Before the proof, notice that the existence of $A^*_1$ as above is equivalent to the condition ${A^*}^{00}_\emptyset \leq \widetilde{G^*}^{000}_\emptyset \cap A^*$. In particular, if ${A^*}^{000}_\emptyset = {A^*}^{00}_\emptyset$, then $A^*_1:={A^*}^{000}_\emptyset$ satisfies the above requirements.\\[3mm]
{\em Proof.} The fact that  $F \leq_B E_{\widetilde{G^*}}$ follows easily from Proposition \ref{E_0 < E_G} applied to $H:=A^*$ and Fact \ref{Borel section} applied to the natural continuous function from  $(\widetilde{G^*}^{00}_\emptyset \cap A^*)/{A^*}^{00}_\emptyset$ onto $(\widetilde{G^*}^{00}_\emptyset \cap A^*)/A^*_1$.

The rest is the proof of the more important reduction $E_{\widetilde{G^*}} \leq_B F$.
Let $M\prec \G^*$ be countable. We would like to pay the reader's attention to the fact that $S_{\widetilde{G^*}^{000}_\emptyset}(M)$ is not a Polish space, but it is an $F_\sigma$ (hence Borel) subset of the Polish space  $S_{\widetilde{G^*}^{00}_\emptyset}(M)$.
We start from the following claim.\\[3mm]
{\bf Claim} There exists a Borel function $\Psi \colon S_{{G^*}^{00}_\emptyset}(M) \to  S_{\widetilde{G^*}^{000}_\emptyset}(M)$ such that $$\Psi(tp(g/M))=tp((a_g,g)/M)$$ for some $a_g \in A^*$. In particular, $tp(a_g/M)$ depends only on $tp(g/M)$.\\[3mm]
{\em Proof.} For $p(y) \in S_{{G^*}^{00}_\emptyset}(M)$ put 
$$[p]= \{ q(x,y) \in   S_{\widetilde{G^*}^{00}_\emptyset}(M) : q \! \upharpoonright \! y =p \},$$ 
a closed subset of  $S_{\widetilde{G^*}^{00}_\emptyset}(M)$. 

By Fact \ref{G^{000}}, 
$$S_{\widetilde{G^*}^{000}_\emptyset}(M)=\bigcup_{i \in \omega} D_i$$ 
for some closed subsets $D_i$ of $S_{\widetilde{G^*}^{00}_\emptyset}(M)$. 
Since for any closed subset $D$ of  $S_{\widetilde{G^*}^{00}_\emptyset}(M)$, $\{ p \in S_{{G^*}^{00}_\emptyset}(M) : [p]\cap D \ne \emptyset \}$ is a closed subset of $S_{{G^*}^{00}_\emptyset}(M)$, we obtain that for each $i \in \omega$
$$D_i':= \{ p \in  S_{{G^*}^{00}_\emptyset}(M) : [p] \cap D_i \ne \emptyset \wedge (\forall j<i)([p] \cap D_j = \emptyset)\}$$
is a Borel subset of $S_{{G^*}^{00}_\emptyset}(M)$.
Moreover, since for every $p \in S_{{G^*}^{00}_\emptyset}(M)$ the intersection $[p] \cap S_{\widetilde{G^*}^{000}_\emptyset}(M)$ is non-empty (which follows from the assumption ${G^*}^{000}_\emptyset = {G^*}^{00}_\emptyset$), we see that
$$S_{{G^*}^{00}_\emptyset}(M) = \bigcup_{i \in \omega} D_i',$$
and, of course, $D_i'$'s are pairwise disjoint.

Consider the Effros Borel space $F(S_{\widetilde{G^*}^{00}_\emptyset}(M))$ (described in Section 1). Define a function $\Phi \colon S_{{G^*}^{00}_\emptyset}(M) \to F(S_{\widetilde{G^*}^{00}_\emptyset}(M))$ by 
$$\Phi(p)=[p] \cap D_i\;\, \mbox{for}\;\,  p \in D_i'.$$
We check that $\Phi$ is a Borel function. Since all $D_i'$'s are Borel, it is enough to show that $\{ p\in  S_{{G^*}^{00}_\emptyset}(M) : [p] \cap D_i \cap U \ne \emptyset\}$ is Borel for every open subset $U$ of $S_{\widetilde{G^*}^{00}_\emptyset}(M)$. But this is true, because $D_i\cap U$ is of type $F_\sigma$.

Fact \ref{function d} yields a Borel function $d \colon F(S_{\widetilde{G^*}^{00}_\emptyset}(M)) \to S_{\widetilde{G^*}^{00}_\emptyset}(M)$ satisfying $d(F) \in F$ for every non-empty $F \in F(S_{\widetilde{G^*}^{00}_\emptyset}(M))$.

Define $\Psi \colon  S_{{G^*}^{00}_\emptyset}(M) \to  S_{\widetilde{G^*}^{000}_\emptyset}(M)$ to be the composition $d \circ \Phi$; in other words, 
$$\Psi(p)=d([p]\cap D_i)\; \mbox{for}\; p\in D_i'.$$ 
$\Psi(p) \in  S_{\widetilde{G^*}^{000}_\emptyset}(M)$, because $D_i \subseteq  S_{\widetilde{G^*}^{000}_\emptyset}(M)$ and $d([p] \cap D_i) \in D_i$.
Since $\Phi$ and $d$ are Borel functions, so is $\Psi$. The fact that $\Psi(p) \in [p]$ guarantees that $\Psi(tp(g/M))=tp((a_g,g)/M)$ for some $a_g \in A^*$. \hfill $\square$\\[3mm]

Define $f \colon S_{\widetilde{G^*}^{00}_\emptyset}(M) \to (\widetilde{G^*}^{00}_\emptyset \cap A^*)/A^*_1$ by $$f(tp((a,g)/M))=(a-a_g) + A^*_1.$$
To see that the image of $f$ is contained in $(\widetilde{G^*}^{00}_\emptyset \cap A^*)/A^*_1$, consider any $(a,g) \in \widetilde{G^*}^{00}_\emptyset$.  Since $(a_g,g) \in \widetilde{G^*}^{00}_\emptyset$, we get $(a-a_g,e)=(a,g)(a_g,g)^{-1}\in \widetilde{G^*}^{00}_\emptyset$, so $a-a_g \in \widetilde{G^*}^{00}_\emptyset \cap A^*$.

To finish the proof of the proposition, we need to check the following three properties.
\begin{enumerate}
\item $f$ is well-defined.
\item $f$ is a Borel function.
\item $tp((a_1,g_1)/M)\; E_{\widetilde{G^*}} \; tp((a_2,g_2)/M) \iff (a_1-a_{g_1})-(a_2-a_{g_2}) \in \widetilde{G^*}^{000}_\emptyset \cap A^*$.
\end{enumerate}
1. Suppose $(a,g) \equiv_M (a_1,g_1)$. Then $a \equiv_M a_1$ and $a_g \equiv_M a_{g_1}$. Therefore, both $a-a_1$ and $a_g-a_{g_1}$ belong to ${A^*}^{000}_\emptyset \subseteq A^*_1$, and so  $(a-a_g) + A^*_1=(a_1-a_{g_1}) + A^*_1$.\\[3mm]
2. We will present $f$ as a composition of three Borel functions $f_1,f_2,f_3$.
The function $f_1\colon S_{\widetilde{G^*}^{00}_\emptyset}(M) \to S_{ A^*}(M) \times S_{{G^*}^{00}_\emptyset}(M)$ is defined by $$f_1(tp((a,g)/M))=(tp(a/M),tp(g/M)).$$ 
The function $f_2 \colon S_{A^*}(M) \times S_{{G^*}^{00}_\emptyset}(M) \to  S_{ A^*}(M) \times  S_{\widetilde{G^*}^{00}_\emptyset}(M)$ is defined by
$$f_2(tp(a/M),tp(g/M))=(tp(a/M),\Psi(tp(g/M)))=(tp(a/M),tp((a_g,g)/M))).$$
The function $f_3\colon S_{ A^*}(M) \times  S_{\widetilde{G^*}^{00}_\emptyset}(M) \to A^*/A^*_1$ 
is defined by
$$f_3((tp(a/M),tp((a_g,g)/M)))=(a-a_g) + A^*_1.$$
The fact that $f_3$ is well-defined follows as in point 1.

It is clear that $f=f_3 \circ f_2 \circ f_1$ and that $f_1$ is continuous. The fact that $f_2$ is Borel follows from the claim. To see that $f_3$ is also Borel (even continuous), consider any closed subset $D$ of  $A^*/A^*_1$. 
It's preimage by $f_3$ equals
$$\{ (p,q) \in S_{A^*}(M) \times S_{\widetilde{G^*}^{00}_\emptyset}(M) : (\exists a \models p)(\exists (b,g) \models q)(a-b \in \rho^{-1}[D])\},$$
where $\rho \colon A^* \to A^*/A^*_1$ 
is the quotient map. Since this set is clearly closed, we get that $f_3$ is continuous.\\[3mm]
3. It was proven in the course of the proof of \cite[Proposition 2.15]{GiKr} that for $a_g$'s chosen so that $(a_g,g) \in \widetilde{G^*}^{000}_\emptyset$ for all $g \in {G^*}^{000}_\emptyset={G^*}^{00}_\emptyset$ (which is true in our case by the claim), the function $$\Phi :\widetilde{G^*}^{00}_\emptyset/\widetilde{G^*}^{000}_\emptyset \to (\widetilde{G^*}^{00}_\emptyset \cap A^*)/(\widetilde{G^*}^{000}_\emptyset \cap A^*)$$ defined by $\Phi((a,g) \cdot \widetilde{G^*}^{000}_\emptyset)= a-a_g +(\widetilde{G^*}^{000}_\emptyset \cap A^*)$ is an isomorphism of groups.

Since  $$tp((a_1,g_1)/M)\; E_{\widetilde{G^*}} \; tp((a_2,g_2)/M) \iff (a_1,g_1)(a_2,g_2)^{-1} \in \widetilde{G^*}^{000}_\emptyset,$$
the implications $(\Rightarrow)$ and $(\Leftarrow)$ in point 3. follow from the proofs in \cite{GiKr} of the facts that $\Phi$ is well-defined and injective, respectively. \hfill $\blacksquare$

\begin{corollary}\label{conjecture 1 true}
Consider the situation described in the first two sentences of Fact \ref{main result from GiKr}. Assume that $A^*_1 \subseteq \widetilde{G^*}^{000}_\emptyset \cap A^*$, $\widetilde{G^*}^{000}_\emptyset \cap A^*$ is not type-definable and ${G^*}^{000}_\emptyset ={G^*}^{00}_\emptyset$. Then $E_{\widetilde{G^*}} \sim_B E_0$.
\end{corollary}
Notice that the assumption $A^*_1 \subseteq \widetilde{G^*}^{000}_B \cap A^*$ holds in various interesting cases, e.g. if $ A^*_1={A^*}^{00}_B={A^*}^{000}_B$ and the action of $G$ on $A$ is trivial; in particular, if $\G=((A,+), M)$, $G$ is definable in the structure $M$ and it acts trivially on $A$, and $A^*_1={A^*}^{0}$.\\[3mm]
{\em Proof.} Since ${A^*}^{00}_\emptyset \leq A^*_1$, by Propositions \ref{x} and \ref{xx}, we get that 
$$E_0 \leq_B E_{\widetilde{G^*}} \leq_B F,$$
where $F$ is the relation on $(\widetilde{G^*}^{00}_\emptyset \cap A^*)/A^*_1$ of lying in the same coset modulo $(\widetilde{G^*}^{000}_\emptyset \cap A^*)/A^*_1$. So, it remains to show that $F \leq_B E_0$.

By Claim 1 in the proof of \cite[Theorem 2.1]{GiKr}, we know that $(A^*_1 + A_0) \times {G^*}^{000}_\emptyset$ is a subgroup of $\widetilde{G^*}$ containing $\widetilde{G^*}^{000}_\emptyset$. Therefore, $(\widetilde{G^*}^{000}_\emptyset \cap A^*)/A^*_1$ is a finitely generated abelian group. 
By a result of Weiss, see \cite[Corollary 1.20, Theorem 1.5]{JaKeLo}, such $F$ is Borel reducible to $E_{0}$. 
\hfill $\blacksquare$\\


Consider the situation from the first two sentences of Fact \ref{main result from GiKr}, and suppose ${G^*}^{000}_\emptyset ={G^*}^{00}_\emptyset$. Assume additionally that $A^*_1 \subseteq \widetilde{G^*}^{000}_\emptyset \cap A^*$ and that ${G^*}^{000}_\emptyset$ is perfect. It was proven in \cite{GiKr} that under this hypothesis, the conclusion $\widetilde{G^*}^{000}_\emptyset \ne \widetilde{G^*}^{00}_\emptyset$ of Fact \ref{main result from GiKr} implies the assumption (i) of this fact iff the following equivalence is true:
%
$$\widetilde{G^*}^{00}_B \cap A^*\subseteq A^*_1 \iff \widetilde{G^*}^{000}_B \cap A^*\subseteq A^*_1.$$
This is strongly related to the following question (which is Question 2.12 in \cite{GiKr}) .

\begin{question}\label{type-definability}
Does there exist data with the properties described  in the first sentence of Fact \ref{main result from GiKr} together with ${G^*}^{000}_\emptyset = {G^*}^{00}_\emptyset$, and such that $\widetilde{G^*}^{000}_\emptyset \cap A^*$ is type-definable but different from  $\widetilde{G^*}^{00}_\emptyset \cap A^*$?
\end{question}

If the answer to the above question is negative, than the previous equivalence is true, so the conclusion of Fact \ref{main result from GiKr} implies the assumption (i) (under the hypothesis described before Question \ref{type-definability}). If the answer is positive, then putting $A^*_1=\widetilde{G^*}^{000}_\emptyset \cap A^*$, we get a situation satisfying all the requirements described in the first two sentences of Fact \ref{main result from GiKr} together with the conclusion of this fact, but the assumption (i) is not satisfied. More interestingly, we have the following observation.

\begin{proposition}\label{conjecture 1 false}
If the answer to Question \ref{type-definability} was positive, then for the resulting group $\widetilde{G^*}$ we would have $E_{\widetilde{G^*}} \sim_B \Delta_{2^{\mathbb N}}$, which yielding a counterexample to  Conjecture \ref{special case of main conjecture} (and so to Conjecture \ref{main conjecture}). 
\end{proposition}
{\em Proof.} Let $A^*_1=\widetilde{G^*}^{000}_\emptyset \cap A^*$. Then the assumptions of Proposition \ref{xx} are satisfied, so $E_{\widetilde{G^*}} \leq_B F$, where $F$ is the relation on $(\widetilde{G^*}^{00}_\emptyset \cap A^*)/A^*_1$ of lying in the same coset modulo $(\widetilde{G^*}^{000}_\emptyset \cap A^*)/A^*_1$. But, in this case, $F$ is just the equality on the compact, Polish group  $(\widetilde{G^*}^{00}_\emptyset \cap A^*)/A^*_1$, so $E_{\widetilde{G^*}}$ is smooth.

The fact that $\Delta_{2^{\mathbb N}} \leq_B E_{\widetilde{G^*}}$ follows from the assumption that $A^*_1 \ne \widetilde{G^*}^{00}_\emptyset \cap A^*$, Proposition \ref{E_G sim E_L} and Fact \ref{Newelski}. \hfill $\blacksquare$\\

Corollary \ref{conjecture 1 true} and Proposition \ref{conjecture 1 false} can by summarized as follows.

\begin{corollary}
(i) Assume the answer to Question \ref{type-definability} is negative. Consider the situation described in the first two sentences of Fact \ref{main result from GiKr}. Assume that $A^*_1 \subseteq \widetilde{G^*}^{000}_\emptyset \cap A^*$ and ${G^*}^{000}_\emptyset ={G^*}^{00}_\emptyset$. In this situation, if $\widetilde{G^*}^{000}_\emptyset \ne \widetilde{G^*}^{00}_\emptyset$, then $E_{\widetilde{G^*}} \sim_B E_0$.\\
(ii)  Assume the answer to Question \ref{type-definability} is positive, and let $\widetilde{G^*}$ be the resulting group. Then for $A^*_1:=\widetilde{G^*}^{000}_\emptyset \cap A^*$ the assumptions of the first two sentences of Fact \ref{main result from GiKr} are satisfied and ${G^*}^{000}_\emptyset ={G^*}^{00}_\emptyset$, but $E_{\widetilde{G^*}} \sim_B \Delta_{2^{\mathbb N}}$. In particular, Conjecture \ref{main conjecture} is false.
\end{corollary}
{\em Proof.} The only thing that needs to be explained is the fact that $\widetilde{G^*}^{000}_\emptyset \cap A^*$ is not type-definable in point (i). Suppose it is type-definable. Since we assume the answer to Question \ref{type-definability} is negative, we conclude that $\widetilde{G^*}^{000}_\emptyset \cap A^*=\widetilde{G^*}^{00}_\emptyset \cap A^*$. This together with \cite[Proposition 2.8]{GiKr} applied to a new $A^*_1$ defined as $\widetilde{G^*}^{000}_\emptyset \cap A^*$ yields $\widetilde{G^*}^{000}_\emptyset =\widetilde{G^*}^{00}_\emptyset$, a contradiction.  \hfill $\blacksquare$\\

%
The next fact is \cite[Corollary 2.14(ii)]{GiKr}, and it says that the answer to Question \ref{type-definability} is negative under the additional assumption that $\widetilde{G^*}^{000}_\emptyset \cap A^*$ is an intersection of definable subgroups of finite index in $A^*$. 
Moreover, \cite[Corollary 2.14(ii)]{GiKr} yields the negative answer in the case when the group $G$ is absolutely connected (in the sense of \cite{Gi}) of finite commutator width.

\begin{fact}
Consider the situation from the first sentence of Fact \ref{main result from GiKr}, and assume that ${G^*}^{000}_\emptyset = {G^*}^{00}_\emptyset$. Then, if $\widetilde{G^*}^{000}_\emptyset \cap A^*$ is an intersection of definable subgroups of finite index in $A^*$, then $\widetilde{G^*}^{000}_\emptyset \cap A^*=\widetilde{G^*}^{00}_\emptyset \cap A^*$
\end{fact}

\section{The Lascar group}

Here we briefly consider the Lascar groups, $Gal_{L}(T)$, $Gal_{0}(T)$, in the light of Borel cardinality issues.  We are working in a monster model $\C$ of a first order theory $T$ in a countable language. In early papers, $Gal_{L}(T)$ has been described as 
$Aut(\C)/Autf_L(\C)$, where $Autf_{L}(\C)$ is the subgroup of $Aut(\C)$ generated by the pointwise stabilizers of elementary (small if one wishes) substructures. Now, $\sigma\in Autf_{L}(\C)$ precisely if $\sigma$ fixes all Lascar strong types over $\emptyset$ of countable tuples (see \cite{Lascar-Pillay}), and one recovers the description in the introduction of $Gal_{L}(T)$ as $Aut(\C)$ considered as acting on all the Lascar strong types. \\

Let $M$, $N$ be countable models (elementary substructures of $\C$), let ${\bar n}$ be some enumeration of $N$. We will abuse notation a bit by defining
$S_{{\bar n}}(M)$ to be the space of complete types  over $M$ which are extensions of $tp({\bar n}/\emptyset)$ (to write $S_{tp({\bar n})}(M)$ would be more accurate).  

Let $\nu: S_{\bar n}(M) \to Gal_L(T)$ be the surjection given by 
$$\nu(tp(f({\bar n})/M))=f\cdot Autf_L(\C)$$
for $f \in Aut(\C)$ (notice that the coset $f \cdot Autf_L(\C)$ does not depend on the choice of $f$ as long as $tp(f({\bar n})/M)$ is fixed).

One can easily check that $\nu(p)=\nu(q)$ if and only if $p E^{M,{\bar n}}_{L} q$, where $E^{M,{\bar n}}_{L}$ is the relation $E^M_L$  from Section 2 restricted to  $S_{\bar n}(M)$. Thus, $Gal_L(T)$ can be identified with the quotient of the Polish space $S_{\bar n}(M)$ by the $K_{\sigma}$ equivalence relation $E^{M,{\bar n}}_L$, which is how we want to view it as a descriptive set-theoretic object.

The next proposition tells us that the resulting ``Borel cardinality" of $Gal_{L}(T)$ does not depend on the choice of $M, N$ or the enumeration of $N$. 

\begin{proposition}\label{independence of M,n}
Let $M_1,N_1,M_2,N_2$ be countable, elementary substructures of $\C$, and $\overline{n_{1}}$, $\overline{n_{2}}$ enumerations of $N_{1}$, $N_{2}$, respectively.  Then $E^{M_1,\overline{n_{1}}}_L \sim_B E^{M_2,\overline{n_{2}}}_L$.
\end{proposition}
{\em Proof.} The independence of the choice of $M$ follows as in Proposition 2.3, so we can assume that $M_1=M_2=:M$. To show the independence of the choice of $N$, we can assume that $\overline{n_1}$ is a  subtuple of $\overline{n_2}$.

Define the function $\Phi : S_{\overline{n_2}}(M) \to S_{\overline{n_1}}(M)$ to be restriction  to the variables corresponding to $N_{1}$.

It is clear that $\Phi$ is continuous and that $p E^{M,\overline{n_2}}_L q$ implies $\Phi(p) E^{M,\overline{n_1}}_L \Phi(q)$. To see the converse, consider any $f,g \in Aut(\C)$ such that $tp(f(\overline{n_1})/M) E^{M,\overline{n_1}}_L tp(g(\overline{n_1})/M)$. Then $f(\overline{n_1}) E_L g(\overline{n_1})$,  so $f^{-1}hg \in Aut(\C/N_1) \subseteq Autf_L(\C)$ for some $h \in Autf_L(\C)$. Since $Autf_L(\C)$ is a normal subgroup of $Aut(\C)$, we get  $gf^{-1} \in Autf_L(\C)$. Hence, $f(\overline{n_2}) E_L g(\overline{n_2})$, so $tp(f(\overline{n_2})/M) E^{M,\overline{n_2}}_L tp(g(\overline{n_2})/M)$. This shows that $E^{M,\overline{n_2}}_L \leq_B E^{M,\overline{n_1}}_L$. The opposite reduction follows by Fact 1.7(ii) (as in the proof of Proposition 2.3). \hfill $\blacksquare$\\

The above proposition allows us to define the {\it Borel cardinality of $Gal_L(T)$} as the Borel cardinality of the relation $E^{M,{\bar n}}_L$ (on the space $S_{\bar n}(M)$).


$Gal_{KP}(T)$ is $Aut(\C)/Autf_{KP}(\C)$, where $Autf_{KP}(\C)$ is the group of automorphisms which fix every $E_{KP}$-class (equivalently fix every ``bounded hyperimaginary''). As before, 
$Gal_{KP}(T)$ is the quotient of  $S_{\bar n}(M)$, but 
now by a closed equivalence relation, giving $Gal_{KP}(T)$ the structure of a compact Hausdorff group. The inclusion 
of $Autf_{L}(\C)$ in $Autf_{KP}(\C)$ induces a canonical surjective homomorphism from $Gal_{L}(T)$ to $Gal_{KP}(T)$, 
whose kernel we call $Gal_{0}(T)$. $T$ is said to be $G$-compact iff $Gal_{0}(T)$ is trivial. \\

$Gal_{0}(T)$ can also be given a well-defined Borel cardinality. For example fix again countable models $M, N$ and let ${\bar n}$ be an enumeration of 
$N$. Define $S_{\bar n}^{KP}(M)$ to be the closed subspace of $S_{\bar n}(M)$ consisting of those $q$ such that some (any) realization  of $q$ is 
$E_{KP}$ equivalent to ${\bar n}$. We let $F^{M,{\bar n}}_{L}$ be the restriction of $E^{M,{\bar n}}_L$ above to $S_{\bar n}^{KP}(M)$. And $Gal_{0}(T)$ 
can be identified with the quotient $S_{\bar n}^{KP}(M)/F^{M,{\bar n}}_{L}$. As in the last proposition, the Borel cardinality is independent of the choice of $M, N$ and the enumeration of $N$. So we have a well-defined ``Borel cardinality of $Gal_{0}(T)$'' which in Conjecture 3 we have conjectured to be non-smooth ($\geq_{B}E_{0}$) whenever $Gal_{0}(T)$ is not the trivial group. 



\begin{remark}
$Gal_0(T) \leq_{B} Gal_L(T)$.
\end{remark}


\begin{remark} (i) In \cite{Ziegler}, Ziegler notes that for $T$ the many sorted theory from Example 3.1 above, $Gal_{L}(T)$ is isomorphic (as a group) to the quotient of $\prod_{n}\mathbb{Z}/n\mathbb{Z}$ by the subgroup of ``bounded sequences", namely what we  called $\Gamma/B_{\Gamma}$ in Proposition 3.11. One can show, by methods of Section 3, that in fact
$Gal_{L}(T) \sim_{B} \Gamma/B_{\Gamma}$ (itself equivalent to $\ell^{\infty}$). In fact, in this case, $Gal_{KP}(T)$ is trivial, so $Gal_{0}(T)$ is Borel equivalent to $\ell^{\infty}$.
\newline
(ii) In the modification Example 3.3 (expanding by covering maps), one can show similarly that $Gal_{KP}(T)$ is trivial, and $Gal_{L}(T) = Gal_{0}(T)$ is both isomorphic to (as a group) and Borel equivalent to $\hat{S}^{1}/\mathbb{R}$  (so of complexity $E_{0}$). 
\end{remark}

We now return to the structures $M$ and $N$ introduced at the end of Section 1. Namely, $M$ is any old structure and $G$ is a $\emptyset$-definable group in $M$, and $N$ is the structure obtained by adding a new sort $X$ and a regular action of $G$ on $X$. $M^{*}$, $N^{*}$, $G^{*}$ are saturated versions. It is well-known that $Aut(N^{*})$ is canonically isomorphic to the semidirect product of $G^{*}$ and $Aut(M^*)$. Here the action of $Aut(M^{*})$ on $G^{*}$ is is the obvious one. 
Fixing a point $x_{0}\in X$, the isomorphism $F$ say between $Aut(N^{*})$ and $G^{*} \rtimes Aut(M^{*})$ takes $f\in Aut(N^{*})$ to  $(g,f|M^{*})$, where $f(x_{0}) = g^{-1}\cdot x_{0}$.  In \cite{GiNe}, it is observed (Proposition 3.3 there) that $F$ induces an isomorphism between $Gal_{L}(Th(N))$ and $(G^{*}/(G^{*})^{000}) \rtimes Gal_{L}(Th(M))$, as well as between $Gal_{KP}(Th(N))$ and  $(G^{*}/(G^{**})^{00}) \rtimes Gal_{KP}(Th(M))$. One deduces an isomorphism between $Gal_{0}(Th(N))$ and $(G^{*})^{00}/(G^{*})^{000} \rtimes Gal_{0}(Th(M))$. The proof of Proposition 2.11 can be suitably modified to obtain {\em Borel equivalences} above in addition to isomorphisms of groups.   
%
%
In particular, we get
\begin{proposition}
Let $M$,$ N$ and $G$ be as above. Then the Borel cardinality of $E_G$ is less than or equal to the Borel cardinality of $Gal_0(Th(N))$. Moreover, if $M$ happens to have constants for an elementary substructure (so that $Gal_{L}(Th(M))$ is trivial), then the Borel cardinality of $E_{G}$ is equal to that of $Gal_{0}(Th(N))$.
\end{proposition}

We return to Remark 5.2 and generalize it somewhat.
For a (possibly infinite) countable tuple $\overline{a}$ of elements of $\C$ consider the closed set $S_{\overline{a}}(M):=\{ tp(\overline{a}'/M): \overline{a}' \equiv \overline{a}\}$ and its closed subset $S_{\overline{a}}^{KP}(M):=\{ tp(\overline{a}'/M): \overline{a}' E_{KP} \overline{a}\}$. Let $E^{M,\overline{a}}_L$ and $F^{M,\overline{a}}_L$ be the restrictions of the relation $E^M_L$ (defined on $S_X(M)$, where $X$ is the sort of $\overline{a}$) to the sets $S_{\overline{a}}(M)$ and $S_{\overline{a}}^{KP}(M)$, respectively. The argument from the proof of Proposition 2.3 shows that the Borel cardinalities of $E^{M,\overline{a}}_L$ and $F^{M,\overline{a}}_L$ do not depend on the choice of $M$. Thus, we will write $E^{\overline{a}}_L$ and $F^{\overline{a}}_L$, having in mind the relations $E^{M,\overline{a}}_L$ and $F^{M,\overline{a}}_L$ for some countable $M \prec \C$.

\begin{remark}
For any countable tuple $\overline{a}$ one has $F^{\overline{a}}_L \leq_B E^{\overline{a}}_L$.
\end{remark}

As any countable tuple $\overline{a}$ from $\C$ is contained a countable $N \prec \C$, the following conjecture seems reasonable.

\begin{conjecture}\label{nice conjecture}
For any countable tuple $\overline{a}$ from $\C$ the Borel cardinality of $E^{\overline{a}}_L$ is less than or equal to the Borel cardinality of $Gal_L(T)$, and the Borel cardinality of $F^{\overline{a}}_L$ is less than or equal to the Borel cardinality of $Gal_0(T)$. 
\end{conjecture} 




One can generalize Conjecture \ref{nice conjecture} in the following way.

\begin{conjecture}
For any countable tuple $\overline{b}$ containing a tuple $\overline{a}$ one has $E^{\overline{a}}_L \leq_B E^{\overline{b}}_L$ and $F^{\overline{a}}_L \leq_B F^{\overline{b}}_L$.
\end{conjecture}

\section{The category of bounded almost hyperdefinable sets}
In this final section, we give some tentative notions of ``definable maps, embeddings, and isomorphisms" between various quotient structures, in both model theory and topology, and make the link with the Borel point of view.\\

Let us first fix a complete (countable) theory $T$, saturated model $\C$ and ``small" set $A$ of parameters. 
By a bounded almost hyperdefinable set, defined over $A$, we mean something of the form $X/E$, where $X$ is a type-definable (even $*$-definable) over $A$ set and $E$ is an $Aut(\C/A)$-invariant equivalence relation on $X$ with boundedly many (i.e. $\leq 2^{\omega + |A|}$) classes.
As we have discussed earlier, if we fix a small model $M_{0}$ containing $A$, then the map taking $a\in X$ to $tp(a/M_{0}) \in S_{X}(M_{0})$ establishes a bijection between $X/E$ and $S_{X}(M_{0})/E'$ for some equivalence relation $E'$ on $S_{X}(M_{0})$. 
We want to give some reasonable definition of a ``definable isomorphism" between two such bounded almost hyperdefinable sets, as well as a ``definable isomorphism" between two quotients of compact Hausdorff spaces by equivalence relations. In the case of compact Polish spaces, this should be a refinement of Borel equivalence.  Likewise we want to give an appropriate definition of a ``definable map" between these kinds of ``spaces".  
The motivation comes in a sense from model theory and geometry rather than set theory.

When we speak simply of a bounded almost hyperdefinable set, we mean something as above which is defined over some small $A$. 

Let us first clarify one of the above remarks. 
\begin{remark} Let $M$ be a small model, and $X$ a type-definable set over $M$. Then a bounded $Aut(\C/M)$-invariant equivalence relation $E$ on $X$ is the {\em same thing} as an equivalence relation $E'$ on $S_{X}(M)$, in the sense that $E$ induces a canonical equivalence relation $E'$ on $S_{X}(M)$, and that conversely any equivalence relation $E'$ on $S_{X}(M)$ yields a canonical bounded $Aut(\C/M)$-invariant equivalence relation on $X$. Moreover, $E$ is type-definable if and only if $E'$ is closed. 
\end{remark}
{\em Proof.}  
The main point is as usual that $E$ is coarser than the equivalence relation on $X$ of having the same type over $M$, so induces a well-defined equivalence relation $E'$ on $S_{X}(M)$.  Conversely, from an equivalence relation $E'$ on $S_{X}(M)$, let $E$ be defined to hold of $(a,b)$ if $E'$ holds of $(tp(a/M), tp(b/M))$. 

For the last remark, note that the function taking $(a,b)$ to $(tp(a/M), tp(b/M))$ takes type-definable (over $M$) sets to closed sets, and the preimage of a closed set is type-definable over $M$. 
\hfill $\blacksquare$\\

\begin{definition} Let $X/E$, $Y/F$ be bounded almost hyperdefinable sets. 
\newline
(i) By a definable map $f$ from $X/E$ to $Y/F$ we mean  a map $f:X/E \to Y/F$which is induced by a type-definable subset $C$ of $X\times Y$ such that $C$ projects onto $X$ and whenever $(x,y)\in C$, $(x',y')\in C$ and $E(x,x')$, then $F(y,y')$.  (To say $f$ is induced by $C$ means $f(a/E) = b/F$ whenever $(a,b)\in C$.) If $f$ is one-one, we say $f$ is a definable embedding of $X/E$ into $Y/F$. 
\newline
(ii) By a definable isomorphism $f$ between $X/E$ and $Y/F$ we mean a bijection $f$ between $X/E$ and $Y/F$ which is induced by
a type-definable $C\subset X\times Y$ which projects onto both $X$ and $Y$ and such that whenever $(x,y)\in C$ and $(x',y')\in C$, then $E(x,x')$ iff $F(y,y')$. 
\end{definition}

\begin{lemma} Suppose $X/E, Y/F$ are bounded almost hyperdefinable sets defined over a small model $M$. Let $f$ be a definable map from $X/E$ to $Y/F$. Then $f$ is induced by a set which is type-definable over $M$.
\end{lemma}
{\em Proof.} Suppose $f$ is induced by $C$, where $C$ is type defined by $\Sigma(x,y,d)$, $d$ some possibly infinite (but small) tuple. Let $q(z) = tp(d/M)$. Let $\Sigma'(x,y)$ be the set of formulas  (equivalent to) $\exists z(q(z)\wedge \Sigma(x,y,z))$. So $\Sigma'(x,y)$ is over $M$, and if $C'$ is the set type-defined by $\Sigma'$, then $C'$ induces $f$. This is easily checked using the fact that if $tp(a/M) = tp(b/M)$, then $E(a,b)$. Likewise for $F$. \hfill $\blacksquare$\\

We can give the same definitions for equivalence relations on arbitrary compact spaces:
\begin{definition}
Let $X,Y$ be compact (Hausdorff) spaces, and $E, F$ equivalence relations on $X,Y$, respectively. By a definable function $f:X/E\to Y/F$ we mean a function induced by a closed subset $C$ of $X\times Y$ which projects onto $X$. Namely, if $(x,y), (x',y')\in C$ and $E(x,x')$, then $F(y,y')$.  
And the obvious things (as in Definition 6.2) for a {\em definable embedding} of $X/E$ in $Y/F$ and a definable isomorphism between $X/E$ and $Y/F$. 
\end{definition}

\begin{remark} (i) In the context of Definition 6.4, if $E$ and $F$ are closed equivalence relations, then $X/E$ and $Y/F$ are themselves compact Hausdorff spaces, and a definable function between $X/E$ and $Y/F$ is the same thing as a continuous function. Likewise, a definable isomorphism is simply a homeomorphism.
\newline
(ii) (In the context of Definition 6.4.) if $f: X/E \to Y/F$ is a definable isomorphism, and $E$ is closed, then so is $F$ (and by (ii) $f$ is a homeomorphism). 
\newline
(iii) (In context of Definition 6.4.) Suppose that $f:X/E \to Y/F$ is an arbitrary map. If $f$ is induced by a continuous function from $X$ to $Y$ (i.e. 
$E(x,x')$ implies $F(g(x), g(x'))$ for some continuous function $g$), then $f$ is a definable map. Also, assuming $X,Y$ to be Polish spaces, if $f$ is a definable map, then $f$ is induced by a Borel map $g$ from $X$ to $Y$.
\newline
(iv) Suppose now that $X, Y$ are type-definable sets (over small $M$), and that $E, F$ are $Aut(\C/M)$-invariant equivalence relations on $X$, $Y$, respectively. Let $E'$, $F'$ be the corresponding equivalence relations on $S_{X}(M), S_{Y}(M)$, respectively. Then a definable function from $X/E$ to $Y/F$ in the sense of Definition 6.2 corresponds precisely to a definable function from $S_{X}(M)/E'$ to $S_{Y}(M)/F'$ in the sense of Definition 6.4.
\end{remark}
{\em Proof.} This is all fairly obvious. For example in the case of (ii), for all $y,y'\in Y$, $F(y,y')$ iff [there are $x,x'\in X$ such that $(x,x')\in C, (y,y')\in C$ and $E(x,x')$].  If $E$ is closed, then (as all spaces in sight are compact) the condition enclosed by [..] is a closed condition on $Y\times Y$, hence $F$ is closed.

For the second part of (iii), if $C \subseteq X\times Y$ is closed and witnesses the definability of $f$, then by 1.7(ii), let $g$ be a Borel section of the (continuous) surjection $C\to X$ and compose it with the projection on the second coordinate to get a Borel map $h:X\to Y$ which clearly induces $f$. 
\hfill $\blacksquare$\\

By 6.5(iii), the notion of a ``definable embedding" from $X/E$ to $Y/F$ ($X,Y$ compact spaces) lies in between that of a continuous reduction and of a Borel reduction. In \cite{Thomas}, 
Simon Thomas gives examples, where $X/E\leq_{B}Y/F$ but there is no continuous reduction. By modifying Thomas's proof, we show that there is no definable embedding in 
these cases. 

\begin{proposition}
There are Borel equivalence relations $E$ and $F$ on the Cantor set $2^{\mathbb N}$ such that $E\leq_B F$ but there 
does not exist a definable embedding from 
$2^{\mathbb N}/E$ to $2^{\mathbb N}/F$. 
\end{proposition}

\noindent {\em Proof.} Recall that an equivalence relation is called {\em countable} if each of its equivalence classes is countable. (Note this is a bit different from model-theoretic parlance where countable would be taken to mean countably many classes.) The equivalence relations 
$E$ and $F$ that make the proposition true will be countable. 

A function $f\colon X\to Y$, $X,Y$ Polish spaces, is called {\em countably continuous} if there are Borel sets 
$B_n\subseteq X$, $n\in {\mathbb N}$, such that $\bigcup_n B_n=X$ and the restriction of $f$ to each $B_n$ is continuous.  
For $C\subseteq X\times Y$ and $x\in X$, let
\[
C_x = \{ y\in Y\colon (x,y)\in C\}.
\]

1. {\em Let $X$ be a Polish space and let $Y$ be a compact metric space. Let  $C\subseteq X\times Y$ be a closed set that projects onto the first coordinate and is such that 
$C_x$ is countable for each $x\in X$. Then 
there exists a countably continuous function $f\colon X\to Y$ with $f(x)\in C_x$, for each $x\in X$.}

To see 1, let $D_n$, $n\in {\mathbb N}$, list closures of sets from some countable open basis of $X$. Let $A_n$ consist of all $x\in X$ such that 
$C_x\cap D_n$ has exactly one point. It is easy to see, using compactness of $Y$ and countability of each $C_x$, that $A_n$ is Borel, in fact, it is $G_\delta$, and that 
$\bigcup_n A_n = X$. Let $B_n = A_n\setminus\bigcup_{i<n} A_i$. For $x\in B_n$, define $f(x)$ to be the unique point in $C_x\cap D_n$. 
Note that this definition gives $f(x)$ for each $x\in X$ and that $f(x)\in C_x$.  
It is straightforward to check, using compactness, that $f$ is continuous on each $B_n$. 

2. {\em Let $E$ and $F$ be Borel equivalence relations on compact metric spaces $X$ and $Y$, respectively. Assume that $F$ is countable. If there exists 
a definable embedding from $X/E$ to $Y/F$, then there is a countably continuous function 
$f\colon X\to Y$ that is a reduction from $E$ to $F$.} 

To see 2, let $C\subseteq X\times Y$ be a definable embedding from $X/E$ to $Y/F$. Then 
$C$ is compact, it projects onto the first coordinate and, since $F$ is countable, $C_x$ is countable for each $x\in X$. The function 
$f$ given by point 1 above applied to this $C$ clearly fulfills the conclusion of point 2. 

Recall that a function $f\colon X\to Y$ is a homomorphism 
from $E$ to $F$ if for $x,y\in X$ with $xEy$ we have $f(x)Ff(y)$. Let $\leq_T$ be the Turing reduction relation and let $\equiv_T$ be the Turing 
equivalence relation among elements of $2^{\mathbb N}$. By a cone we understand a set of elements of $2^{\mathbb N}$ that are $\leq_T$-above a fixed element 
of $2^{\mathbb N}$.  
The following statement is a slight generalization of Theorem 1.1 in \cite{Thomas}; the single change 
consists of assuming that $\theta$ is only countably continuous rather than continuous.

3. {\em Assume $G$ is a countable subgroup of ${\rm Sym}({\mathbb N})$. Let $E_G$ be the orbit equivalence relation 
of the coordinate permuting action of $G$ on $2^{\mathbb N}$. If $\theta\colon 2^{\mathbb N}\to 2^{\mathbb N}$ is 
a countably continuous homomorphism from $\equiv_T$ to $E_G$, then there exists a cone 
that is mapped by $\theta$ to a single equivalence class of $E_G$.}  

We indicate how to make changes in Thomas's proof of \cite[Theorem 1.1]{Thomas} on page 765 of \cite{Thomas}. Keeping the notation from this proof, but using 
only the assumptions from point 3, we 
have a set $C\subseteq 2^{\mathbb N}$, an equivalence relation $E_H$ on $2^{\mathbb N}$ and a countably continuous function $\psi\colon 2^{\mathbb N}\to C$ that is 
a homomorphism from $\equiv_T$ to $E_H\upharpoonright C$. We also have $(E_H\upharpoonright C)\subseteq \equiv_T$. It follows that 
$\psi$ is a homomorphism from $\equiv_T$ to $\equiv_T$. Let $B_n$, $n\in {\mathbb N}$, 
be Borel sets with $\bigcup_nB_n = 2^{\mathbb N}$ and $\psi\upharpoonright B_n$ continuous for each $n$. Since $\bigcup_n B_n = 2^{\mathbb N}$, there is $n_0$ such 
that $B_{n_0}$ is $\leq_T$-cofinal. By Martin's theorem (see \cite[p.164]{Thomas}), there is a compact set $K\subseteq B_{n_0}$ and a cone $D$ such that $K$ 
intersects every $\equiv_T$-equivalence class represented in $D$. 
Since $\psi$ is continuous on the compact set $K$, 
there is a cone $D'$ such that for all $x\in K\cap D'$, we have $\psi(x)\leq_T x$. We claim that this relation holds for all $x\in D\cap D'$. Fix such an $x$. 
Let $y\in K$ be such that $x\equiv_Ty$. Then $\psi(x)\equiv_T\psi(y)$ because $\psi$ is a homomorphism from $\equiv_T$ to $\equiv_T$.   
Since $\psi(y)\leq_Ty$, and as noted $x\equiv_Ty$ and $\psi(x)\equiv_T\psi(y)$, we get $\psi(x)\leq_T x$. Thus, we have that $\psi(x)\leq_T x$ for all $x$ from a cone and 
this is (in addition to Borelness of $\psi$) what is needed to complete 
the argument from \cite{Thomas} proving point 3. 

Now, as in \cite{Thomas}, consider the countable Borel equivalence relations $\equiv_T$ and $\equiv_1$, where $\equiv_1$ is the recursive isomorphism equivalence relation. 
We have $\equiv_T \leq_B \equiv_1$. Following 
the argument from \cite{Thomas}, point 3 implies that there is no countably continuous reduction from $\equiv_T$ to $\equiv_1$, 
which implies, by point 2, that there is no definable embedding from $2^{\mathbb N}/\equiv_T$ to $2^{\mathbb N}/\equiv_1$.
\hfill $\blacksquare$\\

In the proof above, countably continuous functions were used. The dividing line between functions that are countably continuous and those that are not 
has been studied in the past; see for example \cite{So}. We showed above that there is no countably continuous reduction from $\equiv_T$ to $\equiv_1$. On the other hand, as mentioned 
in \cite[p. 762]{Thomas}, the function $2^{\mathbb N}\ni x\to x'\in 2^{\mathbb N}$, where $x'$ is the Turing jump of $x$, is a reduction from $\equiv_T$ to $\equiv_1$. The combination 
of these two facts proves that the Turing jump is yet another example 
of an interesting Borel, in fact, Baire class 1, function that is not countably continuous; for other natural such examples see \cite{So} and references therein. \\

Finally in this section, we want to mix up our categories and speak reasonably coherently of a ``definable map" between a bounded almost hyperdefinable set $X/E$ (in some saturated structure $\C$) and the quotient $Y/F$ of an arbitrary compact (Hausdorff) space $Y$ by an equivalence relation $F$, and also of maps going in the other direction, as well as definable embeddings and isomorphisms.  This can be done by virtue of:

\begin{lemma} Let $X/E$ be a bounded almost hyperdefinable set. Let $M$, $N$ be small models (over which $X$ is defined and such that $E$ is 
invariant under automorphisms which fix pointwise $M$ or $N$). Let $E_{1}$, $E_{2}$ be the corresponding equivalence relations on $S_{X}(M)$, $S_{X}(N)$ induced by $E$ (as in 6.1).
Then $S_{X}(M)/E_{1}$ is definably isomorphic to $S_{X}(N)/E_{2}$. 
\end{lemma}
{\em Proof.} As in the proof of 2.3.
\hfill $\blacksquare$\\


In particular, if we are in a Polish context, then definable isomorphism is a refinement of Borel bi-reducibility.
Many of the Borel bi-reductions proved in this paper can be seen to be (sometimes with a little more work) definable isomorphisms.  For example in both Propositions 3.6 and 3.7, $G/G^{000}$ is definably isomorphic (as a group too) to 
$\widehat{\mathbb{Z}}/\mathbb{Z}$. In the first case, this is witnessed by a continuous (surjective) map between the relevant compact spaces, but not in the second case. 


In the background are also ``noncommutative quotients" in the sense of noncommutative geometry, and we would guess that definable isomorphism 
implies Morita equivalence, which in turn implies Borel equivalence (in a Polish context). 
\\

 
%


\noindent
{\bf Addresses:}\\[3mm]
Krzysztof Krupi\'nski\\
Instytut Matematyczny, Uniwersytet Wroc\l awski,\\
pl. Grunwaldzki 2/4, 50-384 Wroc\l aw, Poland\\[1mm]
{\bf E-mail:} kkrup@math.uni.wroc.pl\\[3mm]
Anand Pillay\\
School of Mathematics,\\
 University of Leeds, Leeds, LS2 9JT, UK\\[1mm]
 {\bf E-mail:} pillay@maths.leeds.ac.uk\\[3mm]
S\l awomir Solecki\\
Department of Mathematics,\\
University of Illinois at Urbana-Champaign,\\ 
1409 W. Green Street, Urbana, IL 61801, USA\\[1mm]
{\bf E-mail:} ssolecki@math.uiuc.edu

\end{document}